\documentclass[11pt]{article}

\usepackage[cp850]{inputenc}
\usepackage{latexsym,graphicx}

\font\tenmath=msbm10 scaled 1200
\font\sevenmath=msbm7 scaled 1200
\font\fivemath=msbm5 scaled 1200

\newfam\mathfam \textfont\mathfam=\tenmath
\scriptfont\mathfam=\sevenmath \scriptscriptfont\mathfam=\fivemath
\def\math{\fam\mathfam}
\def\R{{\math R}}
\def\N{{\math N}}
\def\E{{\math E}}

\def\P{{\math P}}

\def\Q{{\math Q}}

\def \^#1{\if#1i{\accent"5E\i}\else{\accent"5E#1}\fi}

\def \a{\alpha}
\def \b{\beta}
\def \d{\delta}

\def \n{\eta}

\def \cqfd{\quad_\diamondsuit}
\def \ms{\medskip}
\def \ss{\smallskip}
\def \bs{\bigskip}
\def \ni{\noindent}

\newtheorem{Thm}{Theorem}

\newtheorem{Lem}{Lemma}
\newtheorem{Pro}{Proposition}
\newtheorem{Cor}{Corollary}

\author{\sc
{\sc Siegfried Graf}\thanks{Universit\"at Passau, Fakult\"at f\"ur Mathematik und Informatik, D-94030
Passau, Germany. E-mail: {\tt  graf@fmi.uni-passau.de}} \quad
{\sc Harald Luschgy}\thanks{Universit\"at Trier, FB IV-Mathematik, D-54286 Trier, Germany.
E-mail: {\tt luschgy@uni-trier.de}} \quad {\sc  and}
\quad {\sc Gilles Pag\`es} \thanks{Laboratoire de Probabilit\'es et Mod\`eles al\'eatoires, UMR~7599, Universit\'e Paris 6, case 188, 4,
pl. Jussieu, F-75252 Paris Cedex 5, France. E-mail:{\tt  gpa@ccr.jussieu.fr}}
}
\title{\bf Distortion mismatch in the quantization of probability measures}

\begin{document}

%\begin{center}
%{{\bf \Large Quantization of probability distributions under } } \\ $ $ \\
%{{\bf \Large norm-based distortion measures}} \\
%\end{center}

\maketitle
\begin{abstract}
We elucidate the asymptotics of the $L^s$-quantization error induced by a sequence of $L^r$-optimal $n$-quantizers of a
probability distribution $P$ on $\R^d$ when $s>r$. In particular we show that under natural assumptions, the optimal rate is preserved as
long as $s<r+d$ (and for every
$s$ in the case of a compactly supported distribution). We derive some applications of these results to the error bounds for quantization based
quadrature formulae in numerical integration on $\R^d$ and on the Wiener space. 
\end{abstract}

\bigskip
\noindent {\em Key words:  Optimal quantization, Zador Theorem, numerical integration.}

\bigskip
\ni {\em 2000 Mathematics Subject Classification: 60G15, 60G35, 41A25.}
\section{Introduction}
\setcounter{equation}{0}
\setcounter{Assumption}{0}
\setcounter{Theorem}{0}
\setcounter{Proposition}{0}
\setcounter{Corollary}{0}
\setcounter{Lemma}{0}
\setcounter{Definition}{0}
\setcounter{Remark}{0}
Optimal quantization is devoted to the best approximation in $L_{\R^d}^r (\P)$ ($r>0$) of a random vector $X:(\Omega,{\cal A},\P)\to
\R^d$ by random vectors taking finitely many values in $\R^d$ (endowed with a norm $\|\,.\,\|$). When $X\!\in L^r(\P)$, this leads for
every
$n\ge 1$ to the  following
$n$-level
$L^r(\P)$-optimal quantization problem for the random vector $X$ defined by  
\begin{equation}\label{enrX}
e_{n,r}(X):= \inf\left\{\|X-q(X)\|_{_r},\; q:\R^d\to\R^d, \mbox{Borel function},\; {\rm card }(q(\R^d))\le n\right\}
\end{equation}
One shows that the above infimum can be taken over the   the Borel functions 
\[
q:\R^d\to \a:= q(\R^d),\quad \a \subset \R^d,\quad {\rm card} \a \le n
\]
which are some {\em projection following the nearest neighbour rule} on their image $i.e.$
\[
q(x) = \sum_{a\in \a} a \,\mbox{\bf 1}_{V_a(\a)},
\]
$(V_a(\a))_{a\in\a}$ being a Borel partition of $\R^d$ satisfying
\[
V_a(\a) \subset \left\{x\!\in\R^d\,:\, \|x-a\|= \min_{b\in \a} \|x-b\|\right\}.
\]
The set $\a=q(\R^d)$ is (also) called a {\em Voronoi $n$-quantizer} and one denotes 
\[
\widehat X^\a:= q(X).
\]
Then, if $d(x,\a):=\min_{a\in\a}\|x-a\|$ denote the distance of $x$ to the set $\a$, one has  
\[
\|X-\widehat X^\a\|^r_{_r}= \E \,d(X,\a)^r= \int_{\R^d} d(x,\a)^r\P_{_X}(dx)
\]
which shows that $e_{n,r}(X)$ actually only depends on the distribution $P=\P_{_X}$ of $X$ so that 
\[
e_{n,r}(X)= e_{n,r}(P)= \inf_{{\rm card}(\a)\le n} \left(\int d(x,\a)^rdP(x)\right)^{\frac 1r}.
\]
The first two basic results in optimal quantization theory are the following (see~\cite{GRLU1}):

\ss
--  The above infimum is in fact a maximum: there exists for
every $n\ge1$ (at least) one  $L^r(\P)$-optimal $n$-quantizer $\a_n^{*}$ and if ${\rm supp}(P)$ is infinite ${\rm
card}(\a_n^{*})=n$. 

\ss
-- Zador's Theorem: If $X\!\in L^{r+\n}(\P)$ $i.e.$ $\int_{\R^d} \|x\|^{r+\n}dP(x) <+\infty$   for some $\n>0$, then 
\[
\lim_n n^{\frac 1d} e_{n,r}(P) = (Q_r(P))^{\frac 1r} \!\in\R_+. 
\]
A more explicit expression is known for the real constant $Q_r(P)$ (see~(\ref{(2.3)}) below). In particular, $Q_r(P)>0$ if and
only if
$P$ has an absolutely continuous part (with respect to the Lebesgue measure $\lambda_d$ on $\R^d$).  When 
$P$ has an absolutely continuous part,  sequence $(\a_n)_{n\ge1}$ of
$n$-quantizer is   $L^r$-{\em rate optimal} for $P$ if 
\[
\limsup_n n^{\frac 1d}\int_{\R^d} d(x,\a_n)^r\P_{_X}(dx)<+\infty
\]
and is {\em asymptotically $L^r$-optimal} if $\displaystyle  \lim_n \frac{\int_{\R^d} d(x,\a_n)^r\P_{_X}(dx)}{e_{n,r}(P)}=1$.  
 
\medskip
Our aim in this paper is to deeply investigate the (asymptotic) $L^s$-quantization error induced  by a sequence $(\a_n)_{n\ge 1}$ of
$L^r$-optimal
$n$-quantizers.  It follows from the monotony of $s\mapsto \|\,.\|_s$ that $(\a_n)_{n\ge 1}$ remains  an $L^s$-rate optimal
sequence as long as  $s\le r$. As soon as $s>r$ no such straightforward answer is available  (except for the  the uniform
distribution over the unit interval since the sequence $((\frac{2k-1}{2n})_{1\le k\le n})_{n\ge 1}$ is $L^r$-optimal for every $r>0$). 

\medskip
Our main motivation for investigating this problem comes from the recent application of optimal quantization in Numerical Probability
to numerical integration (see~\cite{BAPA}) and the computation of
conditional expectation ($e.g.$ for the pricing of American options,
see~\cite{PAG}). Let us consider for the sake of simplicity the   case of the error bound in
the quantization based quadrature formulas for numerical integration. Let
$F:\R^d\to \R$ be a ${\cal C}^1$ function with a Lipschitz continuous 
differential  $DF$. It follows from a simple Taylor expansion
(see~\cite{PAG}) that  for any random vector $X$ with distribution
$P=\P_{_X}$ quantized by $\widehat X^\a$ ($\a\subset \R^d$)
\[
\left|\E(F(X))-\E(F(\widehat X^\a))- \E(DF(\widehat X^\a).(X-\widehat X^\a))\right|\le [DF]_{\rm Lip} \E|X-\widehat X^\a|^2.
\]
If $\a$ is an $L^2$-optimal (or {\em quadratic}) quantizer then it is {\em stationary} (see~\cite{PAG} or~\cite{GRLU1}) so that 
\[
 \widehat X^\a= \E(X\,|\widehat X^\a) 
\]
which makes the first order term  vanish since
\[
\E(DF(\widehat X^\a).(X-\widehat  X^\a))= \E(DF(\widehat X^\a).\E(X-\widehat X^\a\,|\,\widehat X^\a ))=0.
\]
Finally, if $(\a_n)_{n\ge 1}$ is a sequence of {\em quadratic} optimal $n$-quantizers
\begin{equation}\label{ordre2}
\left|\E(F(X))-\E(F(\widehat X^\a))\right|\le [DF]_{\rm Lip}(e_{n,2}(P))^2\sim [DF]_{\rm Lip}\frac{Q_2(P)}{n^\frac 2d}.
\end{equation}   
Now, if the Hessian $D^2F$ does exist, is $\rho$-H\" older ($\rho\!\in(0,1]$) and computable, the same approach yields  
\begin{equation}\label{ordre3}
\left|\E(F(X))-\E(F(\widehat X^{\a_n})) -\E((X-\widehat X^{\a_n})^*D^2F(\widehat X^{\a_n})(X-\widehat X^{\a_n}))\right|\le [D^2F]_{\rho}\|X-\widehat
X^{\a_n}\|^{2+\rho}_{2+\rho}.
\end{equation} 
Consequently evaluating the asymptotic behaviour of $\|X-\widehat X^{\a_n}\|_{2+\rho}= (\int_{\R^d} d(x,\a_n)^{2+\rho}dP(x))^{\frac{1}{2+\rho}}$ is
ne\-cessary to evaluate to what extend  the quadrature formula in~(\ref{ordre3}) does improve the former one~(\ref{ordre3}). Similar problems occur
when evaluating the error in the first order quantization based scheme designed for the pricing of multi-asset American options or for non-linear
filtering (see~\cite{BAPAPR},~\cite{SEL}). One also meets such mismatch problems  in infinite dimensions  when dealing with (product)
functional quantization on the Wiener space in order to price path-dependent European options (see the example in Section~\ref{Appli}
and~\cite{PAPR}).

\ms
The paper is organized as follows: in Section~\ref{two}, a lower bound for the $L^s(\P)$-quantization rate of convergence of an asymptotically
$L^r$-optimal sequence $(\a_n)_{n\ge 1}$ of $n$-quantizers is established. It shows in particular that for absolutely continuous distributions $P$
with unbounded support, it is hopeless to preserve the quantization rate $n^{-\frac 1d}$   rate in $L^s$ as soon as $s>r+d$. When $s\le r+d$ it
provides a lower bound which can be finite. We conjecture that, when $(\a_n)_{n\ge 1}$ is $L^r$-rate optimal the lower bound is in fact 
the sharp rate. In Section~\ref{Trois3}, we show several natural criterions on the distribution $P$ which ensure that $(\a_n)_{n\ge 1}$
is $L^s$-rate optimal for a given $s\!\in(r,r+d)$ or even for any $s\!\in(r,r+d)$.  Our criterions are applied to many parametrized
families of distributions on
$\R^d$. We  investigate by the same method in Section~\ref{quatre} the critical case $s=r+d$ and
the super-critical case
$s>r+d$. In Section~\ref{cinq} we show that for compactly supported distributions on the real line the lower bound obtained in Section~\ref{two}
does hold as a sharp rate. Finally, in Section~\ref{Appli}  we apply  our results to the evaluation of errors in numerical integration by
quantization based quadrature formulas in finite and infinite dimensions.

\bs
{\sc Notations :}  $\bullet$ $\| \cdot \|$ will denote a  norm on $\R^d$ and $B(x,r)$ will denote the closed  ball
centered at $x$ with radius $r>0$ (with respect to this norm),  $d(x,A)$ will denote the distance between $x\!\in\R^d$
and a subset $A\subset
\R^d$. 

\ss
$\bullet$ $\lambda_d$ will denote the Lebesgue measure on $\R^d$
(equipped with its Borel $\sigma$-field ${\cal B}(\R^d)$). 

\ss $\bullet$ Let $(a_n)_{n\ge 0}$ and $(b_n)_{n\ge 0}$ be two sequences of positive real numbers. The symbol
$a_n \asymp b_n$ is for $a_n =O(b_n)$ and
$b_n=O(a_n)$ whereas the symbol
$a_n\sim b_n$ means $ a_n=b_n +o(b_n)$ as $n\to \infty$.

\ss $\bullet$ $[x]$ is for the integral part of the  real number $x$.

\ss $\bullet$ $f\propto g$ means that the functions $f$ and $g$ are proportional.

\section{The lower estimate}\label{two}
\setcounter{equation}{0}
\setcounter{Assumption}{0}
\setcounter{Theorem}{0}
\setcounter{Proposition}{0}
\setcounter{Corollary}{0}
\setcounter{Lemma}{0}
\setcounter{Definition}{0}
\setcounter{Remark}{0}

In this section we derive a precise lower bound in the $(r,s)$-problem for
non-purely singular probability distributions $P$ and asymptotically 
$L^r$-optimal quantizers. It is expected to be best possible. 

Let   $r \in (0, \infty)$. Let  $P$ be a probability measure on
$(\R^d,{\cal B}(\R^d))$ satisfying  
\begin{equation}
\int_{\R^d} \| x \|^r dP(x) <+ \infty
\end{equation}
and $\mbox{supp} (P)$ is not finite. Then $e_{n,r} (P) \in (0, \infty)$
 for every $n$ and 
$e_{n,r}(P) \rightarrow 0$ as $n \rightarrow \infty$. A sequence 
$(\alpha_n)_{n \geq 1}$ of quantizers is called asymptotically $L^r$-optimal  for $P$
if $\mbox{card} \, \alpha_n \leq n$ for every $n$ and
\begin{equation}
\int_{\R^d} d(x, \alpha_n)^r dP(x) \sim e_{n,r} (P)^r \; \mbox{ as } \; n
\rightarrow \infty .
\end{equation}
Let  $P^a = f.\lambda_d$ denote the absolutely continuous part of $P$
with respect to $\lambda_d$ and assume that 
$\displaystyle \int_{\R^d} \| x \|^{r+\eta} dP(x) < +\infty$ for some
$\eta > 0$. Then by the Zador Theorem  (see~\cite{GRLU1})
\begin{equation}\label{(2.3)}
\lim_{n \to \infty} n^{r/d} e_{n,r}(P)^r = Q_r(P)
\end{equation}
where
\begin{equation}
Q_r (P) := J_{r,d} \left( \int_{\R^d} f^{d/(d+r)} d
\lambda_d\right)^{(d+r)/d}\in [0,\infty)
\end{equation}
and
\[
J_{r,d} := \inf_{n \geq 1} n^{r/d} e_{n,r} (U([0,1]^d ))^r \in (0,
\infty).
\]
($U([0,1]^d)$ denotes the uniform distribution on the hyper-cube
$[0,1]^d$). Note that the above moment assumption implies that
$\int_{\R^d} f^{d/(d+r)} d \lambda_d<+\infty$. 

\ss
Furthermore, for probabilities $P$ on $\R^d$ with $P^a \not= 0$, 
% and
%$\displaystyle \int_{\R^d} f^{d/(d+r)} d \lambda_d < +\infty$, 
the empirical measures associated to an asymptotically $L^r$-optimal 
sequence $(\a_n)_{n\ge 1}$ of $n$-quantizers satisfies (see~\cite{GRLU1}
Theorem~7.5 or~\cite{DEGRLUPA})
\begin{equation}\label{empiric}
\frac{1}{n} \sum_{a\in \a_n} \d_{a}
\stackrel{w}{\longrightarrow}P_r
\end{equation}
where $P_r$ denotes the
$L^r$-point density measure of $P$   defined by 
\begin{equation}
P_r := f_r .\lambda_d \quad\mbox{ with }\quad  f_r := \frac{f^{d/(d+r)} }{
\int f^{d/(d+r)} d
\lambda_d } .
\end{equation}
Note that the limit  $Q_r(P)$ in the Zador Theorem reads
\[
Q_r (P) = J_{r,d} \int f^{-r/d}_r d P^a .
\]
The quantity that  naturally comes out in the $(r,s)$-problem, $r, \,s
\!\in (0,\infty)$, is 
\begin{eqnarray}
Q_{r,s} (P) & := & J_{s,d} \int f^{-s/d}_r dP^a \\ 
& = &   J_{s,d} \left( \int_{\R^d}
f^{d/(d+r)} d
\lambda_d\right )^{s/d}
\int_{\{f>0\}} f^{1-s/(d+r)} d \lambda_d\in (0,+\infty] . \nonumber
\end{eqnarray}
\begin{Thm}\label{thmliminf}
Let $r, s \in (0, \infty)$. Assume $P^a \not= 0$ and $\displaystyle
\int_{\R^d}\| x \|^{r+\eta} dP(x) <+ \infty$  for some $\eta > 0$. Let
$(\alpha_n)_{n \geq 1}$ be an asymptotically $L^r$-optimal 
%$n$-quantizer
sequence of $n$-quantizers for $P$. Then
\begin{equation}\label{liminf}  
\liminf_{n \to \infty} n^{s/d} \int d(x, \alpha_n)^s dP(x) \geq Q_{r,s} (P) .
\end{equation}
\end{Thm}

Prior to the proof, let us provide a  few comments on this lower bound.

\bs
\noindent {\sc Comments.} $\bullet$ The main corollary that can be directly derived from Theorem~\ref{thmliminf} is
that
\[
\int_{\{f>0\}} f^{1-s/(d+r)} d \lambda_d 
=+\infty
\Longrightarrow \lim_{n \to \infty} n^{s/d} \int d(x, \alpha_n)^s
dP(x)=+\infty 
\]
since then  $Q_{r,s}(P) = +\infty$.   

By contraposition,  a necessary condition for an asymptotically $L^r$-optimal sequence  of quantizers $(\alpha_n)$
to achieve the optimal rate $n^{-s/d}$ for the $L^s$-quantization error
is that  $Q_{r,s}(P) < +\infty$.  But, under the  moment assumption  of Theorem 1 one has the following equivalence 
holds true
\begin{equation}\label{CNoptirs}
 Q_{r,s}(P) < +\infty\Longleftrightarrow \int f^{-\frac{s}{d+r}} dP^a= \int_{\{f>0\}} f^{1- \frac{s}{d+r}} d
\lambda_d <+
\infty 
\end{equation}
since $\int_{\R^d}
f^{d/(d+r)} d
\lambda_d<+\infty$ (see~\cite{GRLU1}). 

In turn, for probability measures $P$ satisfying $\lambda_d(f > 0) =+ \infty$ a necessary condition for 
 the right hand side of~(\ref{CNoptirs}) to be satisfied  is that
\begin{equation}\label{sledplusr}
s < d + r.
\end{equation}
Indeed, if $s\ge d+r$, the following chain of inequalities holds true 
\[
\lambda_d(f > 0) = \int \mbox{\bf 1}_{\{f>0\}}f^{-1}dP^a\le \left(\int \left(\mbox{\bf
1}_{\{f>0\}}f^{-1}\right)^{\frac{s}{d+r}}\!\!dP^a\right)^{\frac{d+r}{s}}=
\left(
\int_{\{ f > 0
\} }\ f^{1 - s/(d+r)} d
\lambda_d\right)^{\frac{d+r}{s}} 
\]
where we used that $p\mapsto \|\,.\,\|_{L^p(P^a)}$ is non-decreasing since $P^a(\R^d)\le 1$. 

On the other hand, still  when $s<d+r$, the following criterion holds  for the finiteness of $Q_{r,s} (P)$: 
\begin{equation}\label{sledplusr2}
\left(\exists\, \vartheta>0, \int_{\R^d} \| x \|^{ ds/(d+r-s)+\vartheta}
dP(x) < +\infty\right)\Longrightarrow Q_{r,s} (P) < +\infty. 
\end{equation}  
Set $\rho = 1-\frac{s}{d+r}\!\in(0,1)$ and $u= \frac{ds}{d+r-s}+\vartheta$.
Then~(\ref{sledplusr2}) follows  from the regular  H{\"o}lder
inequality applied with  
$\tilde{p} = \frac 1 \rho = \frac{d+r}{d+r-s}$ and $\tilde{q} = \frac{1}{1-\rho}= \frac{s}{d+r-s}$,
\begin{eqnarray*}
\int_{B(0,1)^c} f^\rho d \lambda_d 
& \leq & \left( \int_{B(0,1)^c} (f(x)^\rho \| x \|^{u \rho} )^{\tilde{p}} d\lambda_d(x) \right)^{1/\tilde{p}} 
\left( \int_{B(0,1)^c} \|x\|^{-u\rho\tilde{q}} d \lambda_d(x) \right)^{1 / \tilde{q}} \\
 & = & \left( \int_{B(0,1)^c} f(x) \|x\|^u d \lambda_d(x) \right)^\rho \left( \int_{B(0,1)^c} 
\| x \|^{-u \rho/(1-\rho)} d \lambda_d(x) \right)^{1-\rho} <+ \infty
\end{eqnarray*}
using the moment assumption in~(\ref{sledplusr2}) and  $u \rho/(1-\rho)= d+ \vartheta\frac{\rho}{1-\rho}>d$.

\ni $\bullet$ It is generally not true in the general setting of
Theorem~\ref{thmliminf} that
$\displaystyle \lim_{n \to \infty} n^{s/d} \int d(x, \alpha_n)^s dP(x) =
Q_{r,s}(P)$  (see Counter-Example 2 in Section~\ref{Trois3}).
However, one may reasonably conjecture that this limiting result holds
true for sequences $(\alpha_n)$ of {\em exactly} $L^r$-optimal
$n$-quantizers.   Our   result in one dimension
for compactly supported distributions (see~section~\ref{cinq}) supports this conjecture.

\medskip 
\ni $\bullet$ In any case, note that~(\ref{liminf}) improves the
obvious lower bound 
$$
\liminf_{n \to \infty} n^{s/d} \int d(x,\alpha_n)^s dP(x) \geq \liminf_{n \to \infty} n^{s/d} e_{n,s}(P)^s
\ge  Q_s(P).
$$
(The right inequality is always true and comes from the proof of the Zador Theorem, see~\cite{GRLU1}). In fact, one even has that, for every $r,\,
s\!\in(0,+\infty)$,  
\[
Q_{r,s}(P) \geq Q_s(P).
\]
Furthermore, this inequality is strict when $r\neq s$ (except if $f$ is  $\lambda_d$-$a.e.$ constant on $\{f>0\}$).
Set $p = (d+s)/s>1$, 
$q = (d+s)/d>$ and   $a = ds / ( d+r)(d+s)$, $b = (d+r-s)d/(d+r)(d+s)$. Then  the H{\"o}lder
inequality yields 
\begin{eqnarray*} 
(Q_{s}(P))^{\frac{d}{d+s}}&=&\int f^{d/(d+s)} d \lambda_d  =  \int_{ \{f > 0\} } f^a f^b d \lambda_d \\
& \leq & \left( \int f^{ap} d \lambda_d\right)^{1/p} \left( \int_{\{ f > 0\} } f^{bq} d
\lambda_d\right)^{1/q} \qquad\mbox{(``$<$" if $f^{ap}$ and $f^{bq}$ are not proportional)} 
\\ & = & \left( \int
f^{d/(d+r)} d \lambda_d\right)^{s/(d+s)} \left( \int_{\{f>0\}} f^{1 - \frac{s}{d+r}} d
\lambda_d\right)^{d/(d+s)}\\ 
%& = & \left( \int f^{d/(d+r)} d \lambda_d\right)^{s/(d+s)} \left( \int
%f^{-s/(d+r)} dP^a\right)^{d/(d+s)} \\
 &=& (Q_{r,s}(P))^{\frac{d}{d+s}}.
\end{eqnarray*}

%
%Consequently, the above bound~(\ref{liminf}) 
%improves the obvious lower
%bound 
%(See \cite{DEGRLUPA} for the latter inequality.)

%
\ni {\bf Proof of Theorem 1.} First note that, the $r$-moment assumption
on
$P$ implies the finiteness of $\int f^{\frac{d}{d+r}}d\lambda_d$. The
existence of at least one asymptotically $L^r$-optimal sequence
$(\a_n)_{n\ge 1}$ follows from the moment assumption on $P$ (as well as
the finiteness of $\int_{\R^d}f^{\frac{d}{d+r}}d\lambda_d$).  For every
integer $m\ge1$, set
\[
f_m := \sum_{k=0}^{m2^m-1} \frac{k}{2^m} \mbox{\bf 1}_{_{E^m_k}}\qquad \mbox{ with }\qquad E^m_k = \left\{\frac{k}{2^m}\le f<\frac{k+1}{2^m}\right\}\cap B(0,m).
\]
The sequence $(f_m)_{m\ge 1} $ is non-decreasing and converges to $f\mbox{\bf 1}_{\{0\le f<+\infty\}}=f$ $\lambda_d$-$a.e.$. 
Let $I_m:=\{k\!\in\{0,\ldots,m2^m-1\}\,:\, \lambda_d(E^m_k)>0\}$. For every $k\!\in I_m$,
there exists a  closed set   
$A^m_k\subset E^m_k$ satisfying
%(\footnote{This follows from the regularity of the measure
%$P$ and the fact that any open set $O$ satisfies $O = \cup_{\d>0}O_\d$ with $O_\d=\{x\in \R^d\,|\, d(x,^cO)>\d\}$,
%all these sets having a negligible boundary except possibly   countably many of them.})
\[
 \lambda_d(E^m_k\!\setminus\! A^m_k)  \le
\frac{1}{m^32^{m}}.
\]  
Let $\varepsilon_m\!\in(0,1]$ be a  positive real number  such that the closed
sets $\widetilde A^m_k:=\{x\!\in\R^d\,:\, d(x,A^m_k)\le
\varepsilon_m\},\, k\!\in I_m$, satisfy 
\[
\int_{\widetilde A^m_k} f^{\frac{d}{d+r}}d\lambda_d\le (1+1/m)
\int_{ A^m_k} f^{\frac{d}{d+r}}d\lambda_d<+\infty.
\]
%On can $e.g.$ choose $E^m_k= \{x\in \R^d\,|\, d(x,E^m_k)<\varepsilon\}$ for 
%a small enough $\varepsilon_m>0$.   
Set 
\[
\widetilde f_m := \sum_{k=0}^{m2^m-1} \frac{k}{2^m} \mbox{\bf 1}_{_{A^m_k}}.
\]
It is clear that  
\[
\{f_m\neq \widetilde f_m\} \subset  \bigcup_{0\le k \le 
m2^m-1}(E^m_k\!\setminus\!  A^m_{k}).
\]
Hence 
\[
\lambda_d (\{f_m\neq \widetilde f_m\})\le \sum_{k=0} ^{m2^m-1}\frac{1}{m^32^{m}} = \frac{1}{m^2} 
\]
so that 
\[
\sum_{m\ge 1} \mbox{\bf 1}_{\{f_m\neq \widetilde f_m\}}<+\infty \qquad \lambda_d\mbox{-} a.e. 
\]
$i.e.$, for $\lambda_d$-$a.e.$ $x$,  $f_m(x)=\widetilde f_m(x)$ for large
enough
$m$ so that   $\widetilde f_m$ converges to $f$ $\lambda_d$-$a.e.$. Finally, as a
result $\widetilde f_m\le f_m\le f$ and $f_m$ converges to $f$
$\lambda_d$-$a.e.$. Then, for every $n\ge 1$,
\begin{eqnarray}
\nonumber n^{\frac sd} \int_{\R^d}(d(x,\a_n))^s  d P(x) 
&\ge& n^{\frac sd} \int_{\R^d}(d(x,\a_n))^s  \,\widetilde f_m(x)d\lambda_d(x)\\
\label{hhtilde}&=&  n^{\frac sd}\sum_{k=0}^{m2^m-1} \frac{k}{2^m} \int_{A^m_k} (d(x,\a_n))^s   d\lambda_d(x).
\end{eqnarray}

Since all the sets $\widetilde A^m_k,\, k=0,\ldots,m2^m-1$ are bounded (as
subsets of $B(0,m+1)$), there exists for every $m\ge1$ and every
$k\!\in\{0,\ldots,m2^m-1\}$ a finite ``firewall"
$\b^m_k\subset\R^d$ (see~\cite{GRLU1} or Lemma~4.3 in~\cite{DEGRLUPA} and note that $A^m_k\subset (\widetilde
A^m_k)_{\varepsilon_m/2}:=\{x\!\in\R^d\,:\, d(x,(\widetilde A^m_k)^c)>\varepsilon_m/2\}$) such that

\[
\forall\, n\ge 1,\quad   \forall\, x \!\in A^m_k,\quad d(x,\a_n\cup \b^m_k)=d(x,(\a_n\cup\b_k^m)\cap \widetilde A^m_k).
\]
Set $\b^m =\cup_{0\le k\le m2^m-1}\b^m_k$. Then, for every $k\!\in\{0,\ldots,m2^m-1\}$, for every $x \!\in A^m_k$, 
\[
d(x,\a_n)\ge d(x,\a_n\cup \b_k^m)=d(x,(\a_n\cup\b_k^m)\cap \widetilde A^m_k)\ge
d(x,(\a_n\cup\b^m)\cap \widetilde A^m_k).
\]
Set temporarily  $n^m_k := \mbox{card}((\a_n\cup \b^m)\cap \widetilde A^m_k)$. First note that  it is clear that 
\[
\frac{n^m_k}{n}\sim \frac{ \mbox{card}(\a_n\cap \widetilde A^m_k)}{n}\quad \mbox{ as } n\to \infty. 
\]
It follows from the asymptotic $L^r$-optimality of the
sequence $(\a_n)$ and the empirical measure theorem (see~(\ref{empiric})) that 
\begin{equation}\label{nmk}
\limsup_n \frac{ \mbox{card}(\a_n\cap \widetilde A^m_k)}{n} \le  \frac{\int_{\widetilde
A^m_k}f^{\frac{d}{d+r}}d\lambda_d}{\int f^{\frac{d}{d+r}}d\lambda_d}
\end{equation}
so that 
\[
\liminf_n \frac{n}{n^m_k} \ge  \frac{\int f^{\frac{d}{d+r}}d\lambda_d}{\int_{\widetilde A^m_k}f^{\frac{d}{d+r}}d\lambda_d}\ge
\frac{m}{m+1} \frac{\int f^{\frac{d}{d+r}}d\lambda_d}{\int_{A^m_k} f^{\frac{d}{d+r}}d\lambda_d}.
\]  

On the other hand, for every $k\!\in I_m$, 
%$n^m_k\to +\infty$ as $n\to +\infty$ since $\lambda_d(E^m_k)>0$. On
%the other hand  
\[
\int_{A^m_k} (d(x,\a_n))^s  d\lambda_d(x)\ge \int_{A^m_k}  d(x,(\a_n\cup\b^m)\cap \widetilde A^m_k) d\lambda_d(x)\ge 
\lambda_d(A^m_k) e^s_{n^m_k,s}(U(A^m_k))
\]
where $U(A^m_k)$ denotes the uniform distribution over $A^m_k$ (note that the inequality is trivial when $\lambda_d(A^m_k)=0$). Then one may apply 
Zador's theorem  which yields, combined with~(\ref{nmk}),
\begin{eqnarray*}
\liminf_n n^{\frac sd} \int_{A^m_k}(d(x,\a_n))^s d\lambda_d(x) &\ge& \lambda_d(A^m_k)\times\liminf_n \left(\frac{n}{n^m_k}\right)^{\frac sd}\times
\lim_{n^m_k} ((n^m_k)^{\frac 1d} e_{_{n^m_k}}(U(A^m_k)))^s\\
&\ge & \lambda_d(A^m_k)\times\left(\frac{m}{m+1}\frac{\int f^{\frac{d}{d+r}}d\lambda_d}{\int_{A^m_k}
f^{\frac{d}{d+r}}d\lambda_d}\right)^{\!\frac sd}J_{s,d}\times (\lambda_d(A^m_k))^{\frac sd}\\
&\ge& J_{s,d}\left(\frac{m}{m+1}\right)^{\frac{s}{d}}\left(\int f^{\frac{d}{d+r}}d\lambda_d\right)^{\frac sd}
\left(\frac{\lambda_d(A^m_k)}{\int_{A^m_k} f^{\frac{d}{d+r}}d\lambda_d}\right)^{\frac{s}{d}}\, \lambda_d(A^m_k)\\
&\ge & J_{s,d}\left(\frac{m}{m+1}\right)^{\frac{s}{d}}\left(\int f^{\frac{d}{d+r}}d\lambda_d\right)^{\frac sd}
\left(\frac{k+1}{2^m}\right)^{-\frac{s}{d+r}}\, \lambda_d(A^m_k)
\end{eqnarray*}
with the convention $\frac{0}{0}=0$.

% Note that the above inequality holds true even if $k\!\notin I_m$. 
Consequently,   using~(\ref{hhtilde}) and the super-additivity of $\liminf$   yield that, for every $m\ge 1$, 
\begin{eqnarray*}
\liminf_n n^{\frac sd} \int_{\R^d} (d(x,\a_n))^s dP(x) &\ge& J_{s,d}
\left(\frac{m}{m+1}\right)^{\frac{s}{d}}\left(\int_{\R^d} f^{\frac{d}{d+r}}d\lambda_d\right)^{\frac sd}
\sum_{k=0}^{m2^m-1}\frac{k}{2^m} \left(\frac{k+1}{2^m}\right)^{-\frac{s}{d+r}}\lambda_d(A^m_k)\\
&=& J_{s,d}\left(\frac{m}{m+1}\right)^{\frac{s}{d}} \left(\int_{\R^d}
f^{\frac{d}{d+r}}d\lambda_d\right)^{\frac{s}{d}}\int_{\{f>0\}} \widetilde f_m(\widetilde
f_m+2^{-m})^{-\frac{s}{d+r}}d\lambda_d.
\end{eqnarray*}
 
Now, by Fatou's Lemma, one concludes by letting $m$ go to infinity that 
\[
\liminf_n n^{\frac sd} \int_{\R^d} (d(x,\a_n))^sdP(x) \ge J_{s,d} \left(\int_{\R^d}
f^{\frac{d}{d+r}}d\lambda_d\right)^{\frac{s}{d}}
\int_{\{f>0\}}  f^{1-\frac{s}{d+r}}d\lambda_d.\qquad 
\cqfd
\]
%
%
%\newpage

\section{The upper estimate}\label{Trois3}
\setcounter{equation}{0}
\setcounter{Assumption}{0}
\setcounter{Theorem}{0}
\setcounter{Proposition}{0}
\setcounter{Corollary}{0}
\setcounter{Lemma}{0}
\setcounter{Definition}{0}
\setcounter{Remark}{0}

Let $r, s \in (0, \infty)$, $s>r$. In this section we investigate the upper bound 
$\int d(x, \alpha_n)^s d P(x) = O(n^{-s/d})$ for $L^r$-optimal $n$-quantizers $\alpha_n$ which is less precise than the lower bound given in 
the previous section. The reason for the restriction to (exactly) 
optimal $n$-quantizers instead of only asymptotically optimal $n$-quantizers will become clear soon. See $e.g.$ the subsequent
example~2. The formal restriction $s>r$  is only motivated by the fact that the above upper bound is trivial when $s\le r$
since the
$L^p(P)$ -norms are  nondecreasing as a function of $p$ so that 
\[
n^{s/d}\int d(x, \alpha_n)^s d P(x) \le \left(n^{r/d}\int d(x, \alpha_n)^r d P(x)\right)^{\frac sr}.
 \]   

For a sequence $(\alpha_n)_{n\ge 1}$ of finite codebooks in $\R^d$ and $b \in (0, \infty)$ we introduce the maximal
function
$\psi_{_b} : \R^d \rightarrow \R_{+} \cup \{ \infty \}$ by
\begin{equation}
\psi_{_b} (x) := \sup_{n \geq 1 } 
\frac{\lambda_d(B(x, b d(x, \alpha_n))) }{ P(B(x, b d(x, \alpha_n)))}
\end{equation}
(with the interpretation $\frac{0}{0} := 0$). Note that $\psi_{_b}$ is Borel-measurable 
and depends on the underlying norm on $\R^d$. The theorem below provides a criterion based on these maximal functions
that ensures the $L^s$-rate optimality of $L^r$-optimal $n$-quantizers. In Corollary~\ref{Cor1} we derive a more
applicable criterion which only involves the distribution $P$. 

\begin{Thm}\label{thm2} Let $r, s \in (0, \infty)$, $s>r$.
Assume $P^a \not= 0$ and $\int \| x \|^{r+ \n} dP(x)  <+ \infty$ for some $\n > 0$. For
every $n \ge 1$, let $\alpha_n$ be an $L^r$-optimal $n$-quantizer for $P$.
Assume that the maximal function associated with the sequence $(\alpha_n)$
satisfies
\begin{equation}\label{(3.2)}
\psi_{_b}^{s/(d+r)} \in L^1 (P) 
\end{equation}
for some $ b \in (0, 1/2)$. Then
\begin{equation}\label{limsup}
\sup_{n } n^{s/d} \!\int \! d(x, \alpha_n)^s dP(x)  <+
\infty .
\end{equation}
\end{Thm}

\ni{\bf Remarks.} $\bullet$ In fact the proof of Theorem~\ref{thm2} provides a bit more information than~(\ref{limsup}): 
it follows from~(\ref{(3.6)}) below and   the
(reverse) Fatou's lemma for sequences of dominated non-negative functions 
%Fatou lemma 
that  there exists a real constant $C_3(b)>0$ such that 
\[
\limsup_{n \to \infty} n^{s/d} \!\int\! d(x, \alpha_n)^s dP(x) \leq \int \limsup_{n \to \infty} n^{s/d} d(x, \alpha_n)^s dP(x)\le  
C_3(b) \int\!\psi_{_b}(x)^{s/(d+r)} dP(x)<+\infty.
\] 

\ni $\bullet$ Condition~(\ref{(3.2)}) relates the upper estimate for $\int \! d(x,
\alpha_n)^s dP(x)$ to lower bounds for 
$P$-probabilities of (small) balls. If $P^a = P$ it comes rather close to the minimal condition 
$\displaystyle \int \!f^{-s/(d+r)} dP <+ \infty$ (see Theorem 1). 
Note that $\psi^{s/(d+r)}_{_{b_0}} \in L^1(P)$ implies 
$\psi^{s/(d+r)}_{_b} \in L^1(P)$ for every $b \geq b_0$.
\\ \\
{\bf Proof.} Let $y \in \R^d$ and set $\delta = \delta_n = d(y, \alpha_n)$. 
For every $x \in B(y, \delta/2)$ and $a\! \in \alpha_n$, we have 
$\| x - a \| \geq \| y - a \| - \| x - y \|  \geq \delta/2$ 
and hence
\[
d(x, \alpha_n) \geq \delta/2 \geq \| x - y \| , \;x \in B(y, \delta/2).
\]
Let $\beta = \beta_n = \alpha_n \cup \{ y \}$. Then
\[
d(x, \beta) = \| x - y \| , \;x \!\in B(y,    \delta/2) .
\]
Consequently, for every $b \in (0, 1/2)$,
\begin{eqnarray*}
e_{n,r} (P)^r - e_{n+1,r} (P)^r & \geq & \int d(x, \alpha_n)^r dP(x) - \int d(x, \beta)^r dP(x) \\
& \geq & \int_{B(y, \delta b)} ( d(x, \alpha_n)^r - d(x, \beta)^r ) dP(x) \\
& = &    \int_{B(y, \delta b)} ( d(x, \alpha_n)^r -\| x - y \|^r) dP(x) \\
& \geq & \int_{B(y, \delta b)} ( (\delta/2)^r - (b \delta)^r ) dP(x) \\
& = & ((1/2)^r - b^r) \delta^r P(B(y, b \delta)) .
\end{eqnarray*}

One derives that
\begin{equation}\label{(3.3)}
d(y, \alpha_n)^r \leq \frac{C(b) }{ P(B(y, b d(y, \alpha_n)) }  (e_{n,r}(P)^r - e_{n+1,r} (P)^r) 
\end{equation}
for every $y \in \R^d, b \in (0,1/2), n \ge 1$, where $C(b) = ((1/2)^r -
b^r)^{-1} $. Note that $e_{n,r}(P)^r - e_{n+1,r}(P)^r > 0$ for every $n
\in \N$ (see \cite{GRLU1}).

Now we estimate the increments $e_{n,r}(P)^r - e_{n+1,r}(P)^r$. (This extends a corresponding estimate in 
\cite{GRALU2} to distributions with possibly unbounded support.)

Set   $e_{n,r} = e_{n,r} (P)$ for convenience. Let 
$\{V_a : a \in \alpha_{n+1} \} $ with $V_a = V_a (\alpha_{n+1})$ be a
Voronoi partition of $\R^d$ with respect to 
$\alpha_{n+1}$.
Then $P(V_a) > 0$ for all $a \in \alpha_{n+1}$ and  $\mbox{card} \;
\alpha_{n+1}\! =\! n + 1$ (see \cite{GRLU1}),
\[
\mbox{card} \left\{ a \in \alpha_{n+1} : \int_{V_a} \| x - a \|^r dP(x)
>  
\frac{4 e_{n+1,r}^r }{n+1} \right\} \leq \frac{n+1}{4}
\]
and
\[
\mbox{card} \left\{ a \in \alpha_{n+1} : P(V_a) > \frac{4}{n+1} \right\} 
\leq \frac{n+1}{4} .
\]
This implies that
\[
 \beta_{n+1} := \left\{ a \in \alpha_{n+1} : \int_{V_a} \| x - a
\|^r dP(x) \leq
\frac{4 e_{n+1,r}^r }{n+1} , P(V_a) \leq \frac{4}{n+1} \right\}
\]
satisfies $\mbox{card} \beta_{n+1} \geq (n+1)/2$. Choose a closed 
hyper-cube $K = [-m,m]^d$ such that 
$P_r(K) > 3/4$. The empirical measure theorem (see~(\ref{empiric}) above
or \cite{DEGRLUPA,GRLU1}  for
details)  implies 
\[
\lim_{k\to \infty} \frac{\mbox{card}( \; \alpha_k \cap K ) }{k} = P_r (K)  
\]
 since $P_r(\partial K)= \lambda_d(\partial K)=0$. We deduce $\mbox{card}(
\alpha_{n+1}
\cap K )\geq 3(n+1)/4$ and hence
$\mbox{card} ( \beta_{n+1} \cap K )  \geq (n+1)/4$ for large enough  $n$.
Since one can find a tessellation of $K$ into $[(n+1)/8] \vee 1$ cubes of diameter less than 
$C_1 n^{-1/d}$, there exist $a_1, a_2 \in \beta_{n+1}$, $a_1 \not= a_2$
such that
\[
\| a_1 - a_2 \|  \leq C_1 n^{-1/d}
\]
for every $n \geq 3$.
Let $\gamma = \alpha_{n+1} \setminus \{ a_1 \}$. Using
\[
d(x, \gamma) \leq \| x - a_2 \|  \leq \| x - a_1 \|  + \| a_1 - a_2 \|  ,
\]
one obtains
\begin{eqnarray*}
e_{n,r}(P)^r - e_{n+1,r}(P)^r & \leq & \int d(x, \gamma)^r dP(x) - \int d(x,
\alpha_{n+1})^r dP(x) \\
 & = & \sum\limits_{a \in \gamma} \int_{V_a} \| x - a \|^r dP(x) + \int_{V_{a_{1}}} 
       d(x, \gamma)^r dP(x) - \sum\limits_{a \in \alpha_{n+1}} \int_{V_a}
\| x - a \|^r dP(x) \\
 & = &  \int_{V_{a_{1}}}     (d (x, \gamma)^r - \| x - a_1 \|^r ) dP(x) \\
 & \leq &  (2^r-1) \int_{V_{a_{1}}} \| x - a_1 \|^r dP(x) +  2^r \| a_1 - a_2 \|^r 
         P(V_{a_{1}} ) \\
  & \leq & \displaystyle{ \frac{4 (2^r-1) e^r_{n+1,r} }{n+1} + \frac{4 \cdot 2^r C^r_1 n^{-r/d} }{n+1} } .
 \end{eqnarray*}

Consequently, using~(\ref{(2.3)}),
\begin{equation}\label{(3.4)}
e_{n,r} (P)^r - e_{n+1,r}(P)^r \leq C_2 n^{-(d+r)/d}
\end{equation}
for every $n \in \N$ and some finite constant $C_2$ independent of $n$.
Combining~(\ref{(3.3)}) and~(\ref{(3.4)}), we get
\begin{eqnarray}
n^{s/d} d(x, \alpha_n)^s & \leq & C_3 (b) \left( \frac{\lambda_d(B(x, b d(x, \alpha_n)}{P(B(x, b d(x, \alpha_n))} 
\right)^{s/(d+r)}\nonumber \\ 
& \leq & C_3(b)\,\psi_{_b}(x)^{s/(d+r)} \label{(3.6)}
\end{eqnarray}
for every $x \in \R^d, n \in \N, b \in (0,1/2)$ and some finite constant $C_3(b)$. The proof is completed by integrating both sides with
respect to $P$.
\hfill{$\Box$}
 
\bs
Next we observe that in case $s < d+r$ a local version of
condition~(\ref{(3.2)}) is always satisfied. 
\begin{Lem}\label{Lem1} Assume $\int \| x \|^r dP(x) <+ \infty$ for some
$r
\in (0,
\infty)$. Let $(\alpha_n)$ be a sequence of finite codebooks in
$\R^d$ satisfying $\int d(x, \alpha_n)^r dP(x) \rightarrow 0$. Then the associated maximal functions $\psi_{_b}$ are  
locally in
$L^p(P)$ for every $p \in (0,1)$ $i.e.$ 
\[
\forall\, M, b \!\in (0, \infty),\qquad \int_{B(0,M)}\!\! \psi^p_{_b} \,dP <+
\infty.
\]
\end{Lem}
{\bf Proof.} Let $M, b \in (0, \infty)$ and set $A = \mbox{supp} (P)$. Then $\max_{x \in B(0,M) \cap A} d(x, \alpha_n) \rightarrow 0$
(see \cite{DEGRLUPA}) and hence
\[
C(M) := \sup_{n \geq 1} \max_{x \in B(0,M) \cap A} d(x, \alpha_n) <+ \infty .
\]
One derives that
\[
B(x, bd(x, \alpha_n)) \subset B(0, b C(M) + M)
\]
for every $x \in B(0, M) \cap A, n \in \N$.
Define the Hardy-Littlewood maximal function $\varphi :  \R^d \rightarrow \R_{+} \cup \{ \infty \}$ 
with respect to the finite measures $\lambda_d( \cdot \cap B( 0, b C(M) + M))$ and $P$ by
\[
\varphi(x) = \varphi_{_{b},M}(x) := \sup_{\rho > 0} \frac{\lambda_d (B(x,\rho) \cap B(0, b C(M) + M)) }{ P(B(x, \rho))} .
\]
Then
\[
\psi_{_b} (x) \leq \varphi(x)
\]
for every $x \in B(0, M) \cap A$. 
From the Besicovitch covering theorem follows the maximal inequality
\[
P(\varphi > \rho) \leq \frac{C_1 \lambda_d(B(0,b C(M) + M)) }{ \rho}\
\]
for every
$\rho > 0$ where the finite constant $C_1$ only depends on $d$ and the underlying norm. (See \cite{GEOM},
Theorem~\ref{thm2}.19. The result in \cite{GEOM} is stated for Euclidean
norms but it obviously extends to arbitrary norms since any two norm on
$\R^d$ are equivalent.) Consequently,

\begin{eqnarray*}
\int_{B(0,M) \cap A} \psi^p_b d P  
& \leq & \int \varphi^p d P = \int^\infty_0 P(\varphi^p > t) dt \leq 1 + \int^\infty_1 P(\varphi^p > t ) dt \\
& \leq & 1 + C_2 \int^\infty_1 t^{-1/p} dt <+ \infty 
\end{eqnarray*}

\ni where $C_2 = C_1 \lambda_d (B(0, b C(M) + M)$. \hfill{$\Box$} \\

In case $s > d+r$ under a mild assumption on the support of $P$ a local version of (3.2) 
holds provided the minimal condition is satisfied locally.  

\begin{Lem}\label{Lem2} Set $A = \mbox{supp}(P)$. Assume $P = P^a = f.
\lambda_d$, $\displaystyle \int\! \| x
\|^r dP(x) <+ \infty$ for some $r \!\in (0, \infty)$, 
$\lambda_d(\cdot\cap A)$ is absolutely continuous with respect to $P$
 and $A$ is a finite union of closed convex sets. Let $(\alpha_n)$ be a
sequence of finite codebooks in $\R^d$ satisfying 
$$\int d(x, \alpha_n)^r
dP(x) \rightarrow 0.
$$
Then for $p \in (1, \infty]$, the associated
maximal functions $\psi_b$ are locally in $L^p(P)$ provided $1/f$  is locally in
$L^p(P)$.   
\end{Lem}  

\ni{\bf Remark.} The absolute  continuity assumption $\lambda_d(\cdot\cap A)\ll P$ does not follow from the absolute 
continuity of $P$: set $P=c \mbox{\bf 1}_{U}.\lambda_d$  where $U= \cup_{n\ge 0} (r_n -2^{-(n+1)},r_n+2^{-(n+1)})$, $\Q= \{r_n,\;
n\ge 0\}$ and $c=1/\lambda_d(U)$. Then ${\rm supp}(P) = \R$ but $\lambda_d \not \ll P$ since $\lambda_d(U^c)=+\infty$ and $P(U^c)=0$.

\bigskip

\ni {\bf Proof.}
Let $M, b \in (0, \infty)$ such that $B(0,M) \cap A \not= \emptyset$. We have 
\[
C = C(M) := \sup_{n \geq 1 } \max_{x \in B(0,M)\cap A} d(x, \alpha_n) <+ \infty
\]
(see \cite{DEGRLUPA}) and hence
\[
\psi_{_b}(x) \leq \sup_{\rho \leq b \,C} \frac{\lambda_d(B(x,\rho)) }{ P(B(x,\rho))}
\]
for every $x \in B(0,M) \cap A$. There exists a constant $C_1 \in (0, \infty)$ such that
\[
\lambda_d (B(x,\rho) \cap B(0,M) \cap A) \geq C_1 \rho^d
\]
for every $x \in B(0,M) \cap A$ and $\rho \leq b C$ since $B(0,M)\cap A$ is finite union of convex sets (see
\cite{GRLU1}). Define the maximal function
$\varphi = \varphi_{f^{-1}} : \R^d \rightarrow \R_{+} \cup \{ \infty \}$ by
\[
\varphi(x) := \sup_{\rho > 0} \frac{ \int_{B(x,\rho)} f^{-1} dP(\cdot \cap B(0,M)) }{ P(B(x,\rho) \cap B(0,M)) }
\]
and note that $d \lambda_d ( \cdot \cap A)/dP = 1/f$. One derives that 
\begin{eqnarray*}
\psi_{_b}(x) & \leq & C_2 \sup_{\rho \leq b \,C} \; \displaystyle{ \frac{\lambda_d(B(x,\rho) \cap B(0,M) \cap A) }{ P(B(x,\rho) \cap
B(0,M)) } } \\ & \leq & C_2 \,\varphi(x)
\end{eqnarray*}
for every $x \in B(0,M) \cap A$ with $C_2 = \lambda_d(B(0,1))/C$. 
By the $L^p(P( \cdot \cap B(0,M))$-boundedness of the maximal operator 
$f^{-1} \mapsto \varphi_{f^{-1}}$ (see \cite{GEOM}, Theorem~\ref{thm2}.19), we obtain
\[
\| \psi_{_b} \mbox{\bf 1}_{B(0,M)} \| _{L^p(P)} \leq C_2 \| \varphi \mbox{\bf 1}_{B(0,M)} 
\|_{L^p(P)} \leq C_3 \| f^{-1} \mbox{\bf 1}_{B(0,M)} \|_{L^p(P)} .
\]
This yields the assertion. \hfill{$\Box$}  

\bs
\ni {\bf Remark.} One can replace the assumption on $A={\rm supp}(P)$ by a local ``peakless" assumption, namely
\[
\forall\, M, \,c>0,\; \inf_{x\in B(0,M)\cap {\rm supp}(P),\,0<\rho\le c} \frac{\lambda_d(B(x,\rho)\cap B(0,M)\cap {\rm
supp}(P))}{\lambda_d(B(x,\rho))} >0 
\]
which can be satisfied by many subsets which are not  finite unions of
closed convex  sets ($e.g.$ if  $A=\, C^c $, $C$ convex set with a non
empty interior).

\bigskip
An immediate consequence of Theorem~\ref{thm2}, Lemma~\ref{Lem1} and
Lemma~\ref{Lem2} concerns distributions $P$ with compact support.

\begin{Cor}\label{Cor1} (Compactly supported distributions) 
Assume that $\mbox{supp} (P)$ is compact and $P^a \not= 0$. 

\ss
\ni $(a)$ For $s\! \in(r,d + r)$, the assertion~(\ref{limsup}) of Theorem~\ref{thm2} holds true.  

\ss
\ni $(b)$ Let $s > d + r$. Assume that $\mbox{supp} (P)$ is a finite union of compact convex sets and  $P = P^a = f.\lambda_d$, 
%$\lambda_d (\cdot \cap A)$ 
%is absolutely continuous with respect to $P$,  
%
\[
\int f^{-s/(d+r)} dP <+ \infty .
\]
Then the assertion of Theorem~\ref{thm2} holds true. If furthermore 
$f \geq \varepsilon > 0 \; \lambda_d \mbox{-}a.s.$ on $\mbox{supp} (P)$,
then the assertion~(\ref{limsup}) of Theorem~\ref{thm2} holds for every $s \in (0, \infty)$.
\end{Cor}

Note that if $s=d+r$, claim~$(b)$ and the monotony of $s\mapsto \|\,.\,\|_s$ hold true if there is a positive real $\d>0$ such that
$\displaystyle \int f^{-\d} dP <+ \infty$. 

For distributions with unbounded support we provide the following condition.  

\begin{Cor}\label{Cor2} (Distributions with unbounded supports)  Let $s\! \in(r,
d + r)$.  Assume  $P^a \not= 0$, $\displaystyle \int \| x \|^{r+ \d}
dP(x) <+ \infty$ for some $\d > 0$  and  
\begin{equation}
\int_{B(0, M)^c} \left( \sup_{t \leq 2 b \|  x \| } 
\frac{\lambda_d(B(x,t)) }{ P(B(x,t))} \right)^{s/(d+r)} d P (x) <+ \infty 
\end{equation}
for some $M, b \in (0, \infty)$. Then the assertion~(\ref{limsup}) of Theorem~\ref{thm2} holds true.   
\end{Cor}

\ni{\bf Proof.} Let $x_0 \in \mbox{supp}(P)$. Then $d(x_0, \alpha_n) \rightarrow 0$ (see
\cite{DEGRLUPA}). For
$\| x \|  > N := \| x_0 \|  + \sup_{n \geq 1} d(x_0, \alpha_n)$, we have 
$d(x, \alpha_n) \leq 2 \| x \| $ for every $n \geq 1$ and thus
\[
\psi_{_b} (x) \leq \sup_{t \leq 2 b \| x \| } 
\frac{\lambda_d(B(x,t)) }{ P(B(x,t))}
\]
for every $x \in B(0, N)^c, b > 0$. The assertion follows from
Theorem~\ref{thm2} and Lemma~\ref{Lem1}. \hfill{$\Box$}  

\bigskip
For distributions with radial tails we obtain a condition which 
is very close to the minimal condition $\int f^{-s/(d+r)} dP <+ \infty$. 

\begin{Cor}(Density with radial tails) \label{Cor3} $(a)$ Assume 
$P =  f. \lambda_d$, $\displaystyle \int \| x \|^{r+ \d} dP(x) <+ \infty$ for some $\d > 0$ and 
$f = h (\|  \cdot \|_0)$ on $B_{\|  \cdot \|_0} (0,N)^c$ with
$h : (R, \infty) \rightarrow \R_{+}$ nonincreasing for some $R \in \R_{+}$ and
$\| \cdot \|_0$ any norm on $\R^d$.
%
%\ss
%\ni $(a)$ 
Let  $s \!\in(r,d + r)$.  If 
\begin{equation}\label{(3.7)}
\int f(cx)^{-s/(d+r)} dP(x) <+ \infty
\end{equation}
for some $c > 1$ (and subsequently for every $c'\!\in(1,c]$), then Assertion~(\ref{limsup})  of Theorem~\ref{thm2} holds true.  

\ss
\ni $(b)$ Assume  $d = 1$ and $s \!\in(1,1 + r)$. Assume $P =  f.\lambda$, $\int |x |^{r + \d} dP(x) <+ \infty$ for some 
$\d > 0$. If $\mbox{supp} (P) \subset [R_0, \infty)$ for some $R_0 \in \R$ and $f_{|(R'_0, \infty)}$ is 
nonincreasing for some $R'_0 \ge R_0$. Assume furthermore~(\ref{(3.7)}) for some $c > 1$. Then the assertion~(\ref{limsup}) of
Theorem~\ref{thm2} holds true.
\end{Cor}

\ni {\bf Proof.} $(a)$ 
We may assume without loss of generality that $\| \cdot \| = \| \cdot \|_0$. 
For $b\! \in (0, 1/2)$, let $M = M(b) = N/(1-2b)$. Then for
$x \in B(0, M)^c, t \leq 2b \| x \| $ and $z \in B(0,t)$, we have
\[
\| x+z \|  \geq \| x\|  - \| z \| \geq
\| x \| (1 - 2b) > N .
\]
Consequently,
\begin{eqnarray*}
P(B(x,t)) & = & \int_{B(0,t)} f (z + x) d \lambda_d (z) = \int_{B(0,t)} h ( \| z+x \| ) d\lambda_d (z) \\
& \geq & h( \| x \|  + t) \lambda_d (B(0,t)) \\
& \geq & h( \| x \|  (1+ 2b)) \lambda_d(B(0,t)) 
\end{eqnarray*}
for every $x \in B(0,M)^c$, $t \leq 2 b \| x \| $. Hence
\[
\sup_{t \leq 2\,b \| x \| } \frac{\lambda_d(B(x,t)) }{ P(B(x,t)) } \leq
\frac{1}{h(\| x \|  (1+2b)) } = \frac{1}{f((1+2b)x)}
\]
for every $x \in B(0,M)^c$ and Corollary~\ref{Cor2} gives the assertion. Item $(b)$ follows similarly.\hfill{$\Box$}

 \bigskip
Although Corollary~\ref{Cor3}  provides an efficient criterion for many families of multi-dimensional distributions since, in practise,
most of them do have radial tails, it is also natural to provide a criterion which does not rely on this assumption. This is the aim
of the next Corollary (which can also treat successfully most usual multi-dimensional distributions). 

\begin{Cor}\label{Cor4}  Assume 
$P =  f. \lambda_d$ and $\displaystyle \int \| x \|^{r+ \d} dP(x) <+ \infty$ for some $\d > 0$. Furthermore, assume 
that ${\rm supp} (P)$ has no peak $i.e.$
\begin{equation}\label{peak}
\kappa_f:=\inf_{x\in {\rm supp} (P),\, \rho>0} \frac{\lambda_d ({\rm supp} (P)\cap B(x,\rho))}{\lambda_d (B(x,\rho))}>0
\end{equation}
and that $f$ satisfies the local growth control assumption: there exist    real numbers
$\varepsilon \ge 0$, $\n\!\in(0,\frac 12)$,
$M,\,C>0$ such that    
\begin{equation}\label{Cond2}
 \forall\, x,\,y \!\in{\rm supp }(P),\; \|x\|\ge M,\;\|y-x\|\le
2\n\, \|x\|\;\Longrightarrow \;f(y) \ge C (f(x))^{1+\varepsilon}.
\end{equation}
Then, for every   $s\!\in(r,\frac{d+r}{1+\varepsilon})$ such that  $\displaystyle \int_{\R^d} 
\frac{dP(x)}{(f(x))^{\frac{s(1+\varepsilon)}{d+r}}}<+\infty$ (if any),  the assertion~(\ref{limsup}) from Theorem~\ref{thm2}
holds true. 
%one has   
%\begin{equation}\label{Conclu2}
% \sup_n n^{\frac sd} \!\int_{\R^d} (d(x,\a_n))^sdP(x) <+\infty.
%\end{equation}
In particular, if~(\ref{Cond2}) holds either for  $\varepsilon=0$ or for every $\varepsilon \!\in(0,\underline
\varepsilon]$ $(\underline \varepsilon>0$), and if
\begin{equation}\label{CondInteg}
\forall\, s\!\in(r,d+r),\qquad  \int_{\R^d}  \frac{dP(x)}{(f(x))^{\frac{s}{d+r}}}=\int_{\{f>0\}} 
 (f(x))^{1-\frac{s}{d+r}}d\lambda_d(x)  <+\infty.
\end{equation}
then~(\ref{limsup}) holds   for every $ s\!\in(r,d+r)$.
\end{Cor}

Note that (if $\lambda_d({\rm supp}(P))=+\infty$)
Assumption~(\ref{peak})   is
$e.g.$ satisfied by any finite intersection of half-spaces, the typical
example being
$\R_+^d$. Furthermore, a careful reading of the proof below shows that it
can be slightly relaxed into: there exists a real $c>0$ such that 
\[
\kappa'_f:=\inf_{x\in {\rm supp} (P)} \left\{\frac{\lambda_d ({\rm supp}
(P)\cap B(x,\rho))}{\lambda_d (B(x,\rho))},\; 0<\rho \le c\,\|x\|
\right\}>0.
\]

\bs \ni {\bf Proof.} Let $x_0\!\in {\rm supp}(P)$. Then, for large enough $n$, $\a_n\cap B(x_0,1)\neq \emptyset$.
Hence 
\[
d(x,\a_n) \le d(x,\a_n\cap B(x_0,1))\le \|x-x_0\|+1\le \|x\|+\|x_0\|+1
\] 
so that $d(x,\a_n) \le 2\,\|x\|$ provided $\|x\|\ge \|x_0\|+1$. We will assume from now on (without loss of
generality) that $M\ge \|x_0\|+1$ in Assumption~(\ref{Cond2}).

Then for every $x\!\in{\rm supp}(P)$, $\|x\|\ge M$ and every $\rho\!\in [0, 2\n \|x\|]$
\begin{eqnarray*}
P(B(x,\rho))&=& \int_{B(0,\rho)}f(x+y) \mbox{\bf 1}_{\{{\rm supp}(P)\}}(x+y)dy \\
&\ge & C (f(x))^{1+\varepsilon}\lambda_d(B(x,\rho)\cap {\rm supp}(P))\\
&\ge& C\kappa_f  (f(x))^{1+\varepsilon}\lambda_d(B(x,\rho))
\end{eqnarray*}
where we used~(\ref{peak}),~(\ref{Cond2}) and $\|x\|\le 2\,\n \|x\|$. Then
\[
\sup_{\rho\le 2\n \|x\|}\frac{\lambda_d(B(x,\rho))}{P(B(x,\rho))}\le \frac{1}{C\kappa_f}
\frac{1}{(f(x))^{1+\varepsilon}}.
\]
Finally one concludes  by Corollary~\ref{Cor2} once noticed that 
\[
\qquad \int_{\|x\|\ge M}\sup_{\rho\le 2\n \|x\|}\left(\frac{\lambda_d(B(x,\rho))}{P(B(x,\rho))}\right)^{\frac{s}{d+r}}dP(x)\le
\frac{1}{(C\kappa_f)^{\frac{s}{d+r}}}\int_{\|x\|\ge M}\frac{dP(x)}{(f(x))^{\frac{s(1+\varepsilon)}{d+r}}}<+\infty.\qquad\qquad  
\hfill{\Box}
\]

\ni {\bf Remark.}  Note that the moment assumption~(\ref{CondInteg}) follows from the more natural moment assumption
\[
\forall\, s\!\in(r,d+r),\qquad\exists\,M\ge 0,\; \exists\, \d>0\quad\mbox{ such that }\quad \int_{\|x\|\ge M}
\|x\|^{\frac{ds}{d+r-s}+\d}f(x)dx<+\infty 
\]
which is of course satisfied by all distributions having polynomial moment at any order. This follows from a standard
application of H\"older inequality (like in the remarks that follow Theorem~\ref{thmliminf}).

\bigskip
At the moment Gaussian distributions are most important for applications of the
$(r,s)$-feature.  

\bigskip
\ni {\bf Example 1} 
$\bullet$ Let 
$P=N(0, \Sigma)$ be the   {\it $d$-dimensional normal distribution} with   positive 
definite covariance matrix $\Sigma$.  Then $P$ satisfies 
condition~(\ref{(3.7)}) from Corollary~\ref{Cor3}$(a)$  with
$h(y) = ((2 \pi)^d\mbox{det} \Sigma)^{-1/2} e^{- y^2/2}$,  
$\| x \|_{0} := \|\Sigma^{-1/2} x \|$ where $\| \cdot \|$ stands for  the
Euclidean norm on $\R^d$ and $M=0$. For $s \!\in(r, d + r)$, choose  $c \!\in (1,\sqrt{(d+r)/s})$.

\ms
\ni $\bullet$  For $s \!\in(d, d + r)$,  the assumptions of Corollary~\ref{Cor3}$(a)$ 
are satisfied by  the {\it hyper-exponential distributions} defined by
\[
f(x) \propto  \exp(-a\|x \|^b), \quad x \in \R^d,\; a,\,  b > 0. 
\]
where $\|\,.\, \|$  denotes any norm on $\R^d$  (possible different of the underlying norm).  Set $h(u) \propto \exp(-au^b)$. Note that
if $d=1$, then the normalizing positive real constant is given by $\kappa_{a,b}= \frac{ba^{1/b}}{\Gamma( 1/b)}$. In fact this holds
true for all distributions of the form
\[
f(x) \propto  \|x\|^c \exp(-a\|x \|^b), \quad x \in \R^d,\; a,\,  b >
0,\; c>-d , 
\]
which $e.g.$ include  the (scalar) {\it double Gamma distributions}
on the real line where
\[
f(x) = \frac{a^c}{2 \Gamma (c)} |x |^{c-1} e^{- a|x |} , \,x \in \R, \; a,\,c > 0.
\]

\ni $\bullet$  As concerns scalar distributions   let us
first mention the   {\it Gamma distributions} 
\[
f(x) = \frac{a^b}{\Gamma(b)} x^{b-1} e^{-ax} \mbox{\bf 1}_{(0, \infty)}(x) , \,x \in \R, \; a,\, b > 0 ,
\]
for which 
the assumptions of Corollary~\ref{Cor3}$(b)$  are satisfied for every $s\!\in(1,1 + r)$. This holds as well for  the {\it Weibull
distributions} 
\[
f(x) = b x^{b-1} \exp(-x)^b \mbox{\bf 1}_{(0, \infty)}(x) , \,x \in \R , \;b > 0, 
\]
the {\it lognormal distributions}  
\[
f(x) = \frac{1}{x \sigma \sqrt{2 \pi} } \exp( - \frac{1}{ 2 \sigma^2} ( \log x - a)^2) \mbox{\bf 1}_{(0, \infty)} (x),
\;x \!\in \R,\;\sigma > 0,   \; a \in \R .
\]
\ms
\ni $\bullet$ The {\it logistic distribution}  
\[
f(x) = \frac{e^x}{(1+e^x)^2} , \;x \in \R .
\]
satisfies Corollary~\ref{Cor3}$(a)$ and so do the symmetric {\it $\rho$-stable distributions} with $\rho \in (0,2)$  provided
$r <\rho$ and $s < \frac{\rho}{(1+\rho)}(1+r)$.

\ms
\ni $\bullet$ For the {\it Pareto distributions} where
\[
f(x) = b x^{-(b+1)} \mbox{\bf 1}_{(1, \infty)}(x)   , \;x \in \R,\; b > 0
\]
the assumptions Corollary~\ref{Cor3}$(b)$ are satisfied provided that $r < b$ and $s < b(1+r)/(b+1)$. 

\ms
\ni $\bullet$ In fact all the above examples of distributions also fulfil the criterion proposed in Corollary~\ref{Cor4} (although it 
turns out to be sometimes slightly more demanding). On the other hand some distributions need naturally to call upon this criterion, especially when
their density is not monotonous. Let $\mu= g.\lambda_1$ be an absolutely continuous probability distribution on $[0,1]$  satisfying
\[
0<C_-\le g\le C_+<+\infty.  
\]
Let $(Y_n)_{n\ge 1}$ denote an i.i.d. sequence of $\mu$-distributed random variables and let $N$ denote a Poisson  random variable with parameter
$\lambda> 1$, independent of  $(Y_n)_{n\ge 1}$. Then set $X=N+Y_N$ and    $P=\P_{_X}$.  One easily checks that 
\[
\P_{_X} =  f.\lambda_{_1}\quad \mbox{ with }\quad f(x) = e^{-\lambda}\frac{\lambda^{[x]}}{[x]!} \,g(x-[x])\mbox{\bf 1}_{\R_+}(x).
\]
 Let $\varepsilon>0$ and $\n\!\in(0,(1/2)\wedge \varepsilon)$ and $x\!\in \R_+$, $y\!\in [x(1-\n), x(1+\n)]$. As soon as $n\ge
[\lambda]+1$,   the sequence $\frac{\lambda^n}{n!}$ is decreasing (to $0$). Hence, if $x\ge (1+\n)[\lambda]+1$
\begin{eqnarray*}
\frac{f(y)}{(f(x))^{1+\varepsilon}} &\ge &\frac{C_-}{(C_+)^{1+\varepsilon}}e^{\lambda \varepsilon} \frac{\lambda^{[y]}}{[y]!} 
\frac{([x]!)^{1+\varepsilon}}{\lambda^{(1+\varepsilon)[x]}}\\
&\ge& \frac{C_-}{(C_+)^{1+\varepsilon}}e^{\lambda \varepsilon}\lambda^{ -1-x( \varepsilon-\n)}\frac{([x]!)^{1+\varepsilon}}{(1+[x(1+\n)])!}.
\end{eqnarray*}
One concludes using the Stirling formula that
\[
\liminf_{x\to +\infty} \inf_{y\in(x(1-\n),x(1+\n))} \frac{f(y)}{(f(x))^{1+\varepsilon}} =+\infty.
\]
On the other hand, Assumption~(\ref{CondInteg}) is fulfilled since $X$ has finite moment at any polynomial order (see the remark following
Corollary~\ref{Cor4}). Note that when the density function $g$ is not non-decreasing $f$ cannot be a non-decreasing function so that Corollary~\ref{Cor3} does
not apply.

\bigskip
The exact optimality assumption made on the sequence $(\alpha_n)$ in Theorem~\ref{thm2} (and
Corollary~\ref{Cor1}) is critical to get the optimal rate $n^{-s/d}$. This
is illustrated by the following counter-example. 

\bigskip
\ni {\bf (Counter-)example 2.}
Let $P = U([0,1])$ and for $n \geq 2$ and $\vartheta \in (0, \infty)$ set
\[
\alpha_n = \alpha_n(\vartheta ) := \left\{ \frac{1}{2 n^\vartheta} \right\} \cup \left\{ \frac{1}{n^\vartheta} +
\left(1 - \frac{1}{n^\vartheta} \right) \frac{2(k-1) -1}{2(n-1)} : k = 2 , \ldots , n \right\} .
\]
Let $r \in (0, \infty)$ and assume $\vartheta > \frac{r}{(1+r)}$. Using a non-Voronoi partition gives the 
upper estimate 
\begin{eqnarray*}
\E |X - \hat{X}^{\alpha_n} |^r & \leq & \int \mbox{\bf 1}_{[0, n^{-\vartheta}]} (x) |x - 
\frac{1}{2n^\vartheta} |^r dx \\
& & +  \sum\limits^{n}_{k=2} 
\int \mbox{\bf 1}_{ [ \frac{1}{n^\vartheta }   + ( 1 - \frac{1}{n^\vartheta })
\frac{k-2}{n-1} , \frac{1}{n^\vartheta}+ (1 - \frac{1}{n^\vartheta}) \frac{k-1}{n-1}] } (x)
\left|x - \left( \frac{1}{n^\vartheta } + ( 1 - \frac{1}{n^\vartheta } ) \frac{2(k-1)-1}{2(n-1)} 
\right)\right|^r dx \\
& = & \frac{2}{r+1} ( \frac{1}{2 n^\vartheta } )^{r+1} + ( 1 - \frac{1}{n^\vartheta } )^{r+1}
\frac{1}{(r+1) 2^r (n-1)^r} .
\end{eqnarray*}
Hence
\begin{eqnarray*}
\limsup_{n \to \infty} n^r \E |X - \hat{X}^{\alpha_n} |^r 
& \leq & \limsup_{n \to \infty} \left( \frac{1}{2^r(r+1)} \frac{1}{n^{(r+1) \vartheta -r} } +
         \frac{1}{(r+1) 2^r } \right) \\
& = & \frac{1}{2^r (r+1)} = J_{r,1} = Q_r (P) 
\end{eqnarray*}
so that in fact
\[
\lim_{n \to \infty} n^r \E |X - \hat{X}^{\alpha_n} |^r = Q_r(P)
\]
It follows that the sequence $(\alpha_n(\vartheta ))_n$ is an asymptotically $L^r$-optimal 
$n$-quantizer for every $\vartheta \in (r/(r+1), \infty)$.
Now, let $s > r$ and $\vartheta \in (r/(r+1), s/(s+1))$. Then
\[
\E |X - \hat{X}^{\alpha_n} |^s 
\geq  \int^{1/2 n^\vartheta }_{0} |x - \frac{1}{2 n^\vartheta } |^s dx 
=  \frac{1}{2^{s+1} (s+1)} n^{- \vartheta (s+1)}
\]
so that 
\[
n^s\,\E\,|X - \hat{X}^{\alpha_n} |^s \geq \frac{1}{2^{s+1} (s+1)} n^{s-\vartheta (s+1)}
\]
Consequently,
\[
\lim_{n \to \infty} n^s \,\E\,|X - \hat{X}^{\alpha_n} |^s = \infty .
\]

We finally comment on what can be said without condition~(\ref{(3.2)}).
\\ \\
{\bf Remark} In the situation of Theorem~\ref{thm2} but without assuming~(\ref{(3.2)}) let $P = P^a$ and $s < d + r$. 
One can deduce from~(\ref{(3.6)}) and differentiation of measures that 
\[
\limsup_{n \to \infty} n^{1/d} d( \,\cdot\,, \alpha_n) \leq C_3(b)^{1/s} f^{-1/(d+r)} <+ \infty \;
\qquad P\mbox{-}\mbox{a.s.}
\]
This improves considerably for absolutely continuous distributions and exactly $L^r$-optimal quantizers and
a.s. result in \cite{DEGRLUPA}. Furthermore, by Lemma 1 and Fatou's lemma
(also without~(\ref{(3.2)})),
\begin{eqnarray*}
\sup_{M > 0} \limsup_{n \to \infty} n^{s/d} \int_{B(0,M)}  d(x, \alpha_n)^s dP(x) 
& \leq & \int \limsup_{n \to \infty} n^{s/d} d(x, \alpha_n)^s dP(x) \\
& \leq & C_3 (b) \int f^{-s/(d+r)} d P 
\end{eqnarray*}
so that the minimal condition $\int f^{-s/(d+r)} dP <+ \infty$ already gives the second inequality of the remark following Theorem~\ref{thm2}. 
In particular, by H{\"o}lder's inequality, the moment condition for $P$ in Theorem~\ref{thm2} implies the minimal
condition for $s \leq r$. The first inequality in this remark always holds for $s = r$ since by Theorem 4.5 in
\cite{DEGRLUPA} and (2.3),
\begin{eqnarray*}
\sup_{M > 0} \lim_{n \to \infty} n^{r/d} \int_{B(0,M)} d(x, \alpha_n)^r dP(x)
& = & \sup_{M > 0} Q_r (P) P_r(B(0,M)) = Q_r(P) \\
& = & \lim_{n \to \infty} n^{r/d} \int d(x, \alpha_n)^r dP(x) .
\end{eqnarray*}
\section{The critical and super-critical cases}\label{quatre}
In this section, we investigate the upper bound for the rate of convergence when $s\ge d+r$. 
\begin{Pro}\label{Pro1} Let $r\!\in (0, \infty)$.
Assume $P^a \not= 0$ and $\displaystyle \int \| x \|^{r+ \d} dP(x)  <+ \infty$ for some $\d> 0$. For
every $n \ge 1$, let $\alpha_n$ be an $L^r$-optimal $n$-quantizer for $P$.  
\ss
\ni 

\ss
\ni $(a)$ {\sc The critical case $s=d+r$:}  Assume there is a real $\underline \vartheta> 0$ such that for every $\vartheta \!\in (0,
\underline \vartheta\wedge s)$, there exists   a real number
$M\ge 0$ such that 
\begin{equation}\label{superpsi0}
\int_{B(0;M)^c} \left(\sup_{r\le 2b\|x\|} \frac{\lambda_d(B(x,r))}{P(B(x,r))}\right)^{1-\frac{\vartheta}{d+r}}\|x\|^\vartheta dP(x)
<+\infty.
\end{equation}
Then, 
\[
\forall\, \varepsilon\!\in(0,1+ r/d),\qquad \limsup_n n^{1+\frac rd
-\varepsilon}\!\int_{\R^d} (d(x,\a_n))^{d+r}dP(x) <+\infty.
\]
 
\medskip
\ni $(b)$ {\sc The super-critical case $s>d+r$:} Assume there is a real $\vartheta\!\in(s-(d+r), s)$ and  a real number $M\ge 0$ such
that 
\begin{equation}\label{superpsi}
\int_{B(0;M)^c} \left(\sup_{r\le 2b\|x\|} \frac{\lambda_d(B(x,r))}{P(B(x,r))}\right)^{\frac{s-\vartheta}{d+r}}\|x\|^\vartheta dP(x)
<+\infty.
\end{equation}
Then, 
\begin{equation}\label{Conclu3}
 \sup_n n^{\frac{s-\vartheta}{d}} \!\int_{\R^d} (d(x,\a_n))^sdP(x) <+\infty.
\end{equation}
\end{Pro}

\ni{\bf Proof.}  $(b)$ It follows from~(\ref{(3.6)}) that, for every $n\ge 1$, $b\!\in(0,1/2)$, 
\[
\forall\, x\!\in\R^d,\qquad
(d(x,\a_n))^{d+r}\le C_b \psi_{_b}(x)n^{-(1+\frac{r}{d})}.
\]
On the other hand, $x_0\!\in{\rm supp}(P)$ being fixed, for large enough $n$, $\a_n\cap B(x_0,1)\neq \emptyset$ so that  
\[
\forall\, x\!\in\R^d,\qquad d(x,\a_n) \le (1+|x-x_0|)\mbox{\bf 1}_{\{|x-x_0|\ge 1\}} + 2\mbox{\bf 1}_{\{|x-x_0|<1\}} \le
C_0(|x|\vee 1).
\]
Let $\psi_{_b}$ be the maximal function associated to $(\a_n)_{n\ge 1}$ and $b\!\in(0,1/2)$.  Consequently
\begin{eqnarray*}
\forall\, x\!\in\R^d,\qquad (d(x,\a_n))^{s}&\le & (d(x,\a_n))^{s-\vartheta}(C_0(|x|\vee 1))^\vartheta\\
&\le & C'_{b} n^{\frac{\vartheta-s}{d}}  (\psi_{_b}(x))^{\frac{s-\vartheta}{d+r}}(\|x\|\vee1)^{\vartheta}.
\end{eqnarray*}

It follows from Lemma~\ref{Lem1} that, for every real number $M\ge 0$ and $b>0$, 
\[
\int_{B(0;M)} \left(\psi_{_b}(x)\right)^{\frac{s-\vartheta}{d+r}}\|x\|^\vartheta dP(x)
<+\infty.
\]
Combined with  Assumption~(\ref{superpsi}), this clearly  implies $ \psi_{_b}^{\frac{s-\vartheta}{d+r}}\!\in L^1(P)$. 

\ss
\ni $(a)$ The above computations are still valid when $s=d+r$. So one concludes by considering $\vartheta$ arbitrarily
close to $0$ which is precisely made possible by Assumption~(\ref{superpsi0}). $\cqfd$

\bs
From this Proposition one easily derives some corollaries in the formerly mentioned settings. We give them in details
for the super-critical case. The adaptation to the critical case is straightforward.

\begin{Cor}\label{Cor5}Assume the global assumption on $P$ in Proposition~\ref{Pro1} holds and that  $P= f.\lambda_d$. Let $\vartheta
\!\in (s-(d+r),s)$. 

\ss
\ni $(a)$ If  $f$ is radial outside a compact subset of $\R^d$ (as defined  in Corollary~\ref{Cor1}) and if 
\[
\int f(cx)^{-\frac{\vartheta-s}{d+r}}\|x\|^{\vartheta}dP(x)<+\infty \mbox{ for some real $c>1$}
\]
then Assumption~(\ref{superpsi}) is fulfilled. 

\ss
\ni $(b)$ If $f$ satisfies Assumption~(\ref{peak}) and Assumption~(\ref{Cond2}) for some parameter $\varepsilon\ge 0$
and if 
\[
\int f(x) ^{-\frac{(s-\vartheta)(1+\varepsilon)}{d+r}}\|x\|^\vartheta dP(x)<+\infty
\]
then Assumption~(\ref{superpsi}) is fulfilled.
\end{Cor}

The above assumptions are clearly fulfilled by the normal distributions (and in fact  most  distributions  
mentioned in Example~1).  Numerical experiments seem to suggest that the critical rate $n^{1+\frac rd}$ cannot be improved for $P$ with unbounded
support, that is,
\[
\forall\, s>0,\qquad \lim_n n^{1+\frac rd}\int_{\R^d} d(x,\a_n)^sdP(x)=+\infty.
\]

\section{A sharp rate for   distributions on compact intervals}\label{cinq}
\setcounter{equation}{0}
\setcounter{Assumption}{0}
\setcounter{Theorem}{0}
\setcounter{Proposition}{0}
\setcounter{Corollary}{0}
\setcounter{Lemma}{0}
\setcounter{Definition}{0}
\setcounter{Remark}{0}

\medskip

\begin{Pro} Let $P=f.\lambda_{_1}$ where $f:[a,b]\to \R_+$ is a
Lipschitz continuous probability density function, bounded away from $0$ on
$[a,b]$. Let $r>0$. Let  $(\a_n)_{n\ge
1}$ be a sequence of  $L^r$-stationary and 
 asymptotically $L^r$-optimal $[a,b]$-valued $n$-quantizers. Then  for every
$s>0$,
\[
\lim_n \left(n^s\int_a^b \min_{a\in \a_n} |x-a|^s dP(x) \right) =Q_{r,s} (P).
%J_{s,1}\left(\int_a^b f^{\frac{1}{r+1}}d\lambda_{_1}\right)^s
%\int_a^bf^{1-\frac{s}{r+1}}d\lambda_{_1}. 
\]
\end{Pro}

Note that we do not assume {\em a priori} that the $n$-quantizers are $L^r$-optimal but only stationary and
asymptotically $L^r$-optimal. 

\bs
\ni {\bf Proof.} It is possible in this scalar setting to number the
elements of a $n$-quantizer with respect to the natural order on the real
line. Furthermore, one may assume for large enough $n$ that $\a_n\subset
(a,b)$ and ${\rm card} \a_n =n$. Namely, we set
\[
\a_n=\{\a_n^1,\ldots,\a_n^n\},\qquad a<\a_n^1<\a_n^2<\cdots<\a_n^{n-1}<\a_n^n<b.
\]
We also set $\a_n^0= a$ and $\a_n^{n+1}=b$ for convenience. For
every $n\ge 1$,  set for every $x\!\in[a,b]$ 
\begin{equation} 
 \varphi_{r,n}(x)=
 n\sum_{k=1}^n \Delta \a^k_n\mbox{\bf 1}_{V_k(\a_n)}(x) 
\end{equation}
where $V_k(\a_n)$ denotes the Voronoi cell of $\a^k_n$ and $\Delta
\a^k_n= \a_n^k- \a_n^{k-1}$, $k\!\in\{1,\ldots, n+1\}$. We know
from~\cite{DEFOPA2} (see  the proof of Theorem~4$(a)$, p.1101) that $h$ being Lipschitz continuous and bounded away from $0$ on
$[a,b]$, there is a positive real constant $C_{a,b,f}$ such that 
\begin{equation}\label{maxmin}
  \max_{1\le k\le n+1} \Delta \a^k_n \le C_{a,b,f} \min_{1\le k\le n+1} \Delta \a^k_n.
\end{equation}
A proof of this last inequality can be found in the appendix. Set $\displaystyle  P_{r,n} = \frac 1n \sum_{k=1}^n \d_{\a^k_n}$ the
empirical measure associated to the
$n$-quantizer
$\a_n$. By the empirical measure theorem (recalled in~(\ref{empiric})), it follows that for large enough
$n$,
\[
 \forall\, a',b'\!\in[a,b],\quad \frac{{\rm card}\{i\,:\, \a^i_n \!\in [a',b']\}}{n} \longrightarrow  
P_r([a',b']) : =  c_{f,\frac{1}{r+1}} 
\int_{a'}^{b'} f^{\frac{1}{r+1}}d\lambda_{_1} 
\]
%\[
%\max_{1\le k\le n+1} \Delta \a^k_n = O(\frac 1n)  \qquad \mbox{ and }\qquad 
%\varphi_n
%\stackrel{U}{\longrightarrow}  c_{f,\frac{1}{r+1}}f^{\frac{r}{r+1}} \quad
%\mbox{ as } \quad n\to +\infty,
%\]
uniformly with respect to $(a',b')\!\in[a,b]^2$ (where $c_{g,\vartheta} = (\int_a^b g^\vartheta d\lambda_1)^{-1}$).
Combining this convergence with~(\ref{maxmin}) and using that $f$ is bounded away from $0$ on $[a,b]$ implies that
there is a real constant $C'_{a,b,f}$ such that, for every $n\ge 1$ and every $ a',\,b'\!\in[a,b]$  
\[
\frac{b'-a'}{C'_{a,b,f}n} \le\min_{1\le k\le n+1, \a_n^k\in[a',b']}\Delta \a_n^k\le   \max_{1\le k\le n+1,
\a^k_n\in[a',b']} \Delta
\a^k_n
\le\frac{C'_{a,b,f}}{n}(b'-a').
\]
A straightforward  application of   the Arzela-Ascoli theorem yields that the sequence $(\varphi_{r,n})_{n\ge 1}$ is
relatively compact for the topology of uniform convergence  and that all its  limiting functions are continuous. Let
$\varphi_{r}$ denote such a limiting function. On the one hand, for every $ a',\,b'\!\in[a,b]$ ,  
\[
\int_{a'}^{b'} \varphi_{r,n} \,d P_{r,n} = \sum_{k=1}^n \Delta\a^k_n\mbox{\bf 1}_{\{\a^k_n
\in[a',b']\}}\longrightarrow b'-a'\quad \mbox{ as } n\to \infty.
\] 
$i.e.$ $\varphi_{r,n}. P_{r,n}$  converges weakly to the Lebesgue measure $\lambda_{_1}$ over $[a,b]$. On the
other hand  $\varphi_{r,n}. P_{r,n}$ converges  weakly to $\varphi_{r}.  P_r =
c_{f,\frac{1}{r+1}}  \varphi_{r} \,f^{\frac{1}{r+1}}.\lambda_{_1}$ since $\varphi_{r}$ is continuous.
Consequently, $\varphi_{r}$ is entirely and uniquely determined by this equation which implies that
\[
\varphi_{r,n} \stackrel{U}{\longrightarrow}  \varphi_{r} :=  \left(\int_a^b \!\!f^{\frac{1}{r+1}} d\lambda_1\right)\,
f^{-\frac{1}{r+1}}\quad \mbox{ as } \quad n\to \infty.
\]
%and $U$ stands for the uniform convergence on
%the compact interval $[a,b]$. Now, for every
Now for every $k=2,\ldots,n-1$, 
%\begin{eqnarray*}
%|\mu(V_k(\a_n))- f(\a^k_n)\times \frac 12(\a_n^{k+1}-\a_n^{k-1})| &\le& [f]_{_{\rm
%Lip}}\int_{\frac{\a^k_n+\a^{k-1}_n}{2}}^{\frac{\a^{k+1}_n+\a^{k}_n}{2}} (x-\a^k_n)^2dx \\
%&=& \frac 13[f]_{_{\rm Lip}}\left((\Delta \a^{k+1}_n )^3+ (\Delta \a^{k}_n )^3\right)\\
%&=& O\left(\frac{1}{n^2}\right).
%\end{eqnarray*}
\begin{eqnarray*}
n^s\left|\int_{V_k(\a^k_n)}\!\!\!|\a^k_n-x|^s f(x)dx-f(\a^k_n)\int_{V_k(\a^k_n)}\!\!\!|\a^k_n-x|^sdx
\right|&\le& n^s [f]_{_{\rm Lip}}\int_{V_k(\a^k_n)}\!\!\!|\a^k_n-x|^{s+1}dx\\
&=& [f]_{_{\rm
Lip}}\frac{n^s}{s+2}\left(\!\left(\frac{\Delta\a^{k+1}_n}{2}\right)^{s+2}\!\!+\left(\frac{\Delta\a^{k}_n}{2}\right)^{s+2}\right)\\
&=&\frac{ [f]_{_{\rm
Lip}}}{(s+2)2^{s+2}n^2}\left(\!\left(\varphi_{r,n}(\a^{k+1}_n) \right)^{s+2}\!\!\!\!+\left(\varphi_{r,n}(\a^{k}_n)
\right)^{s+2} \right)
\end{eqnarray*}
and 
\begin{eqnarray*}
n^s f(\a^k_n)\int_{V_k(\a^k_n)}|\a^k_n-x|^s d x  & =& \frac{f(\a^k_n)}{(s+2)2^{s+1}n}
\left((\varphi_{r,n}(\a^{k+1}_n) )^{s+1}+(\varphi_{r,n}(\a^{k}_n))^{s+1}\right)\\
&=& \frac{J_{1,s} }{2n} \left( f(\a^k_n)\left(\varphi_{r,n}(\a^{k}_n) \right)^{s+1}+ f(\a^{k+1}_n)\left(\varphi_{r,n}(\a^{k+1}_n) \right)^{s+1} 
\right) + O\left( \frac{1}{n^2} \right)
\end{eqnarray*}
since 
$$
|f(\a^{k+1}_n)-f(\a^k_n)|(\varphi_{r,n}(\a^{k+1}_n))^{s+1}\le [f]_{_{\rm Lip}}\Delta
\a^{k+1}_n(\varphi_{r,n}(\a^{k+1}_n))^{s+1}= \frac{[f]_{_{\rm Lip}}}{n}(\varphi_{r,n}(\a^{k+1}_n))^{s+2}\le C/n.
$$ 

 For $k=1$ and $n$, the above equalities are no longer true due to edge effects but both induced errors are $O(1/n)$. Consequently 
\begin{eqnarray*}
\left|n^s \,\E|X-\widehat X^{\a_n}|^s \right.\!\!&\!\!-\!\!&\!\!\left. \frac{J_{1,s}}{2}\frac{1}{n} \sum_{k=1}^n\left(
f(\a^k_n)\left(\varphi_{r,n}(\a^{k}_n)
\right)^{s+1}+ f(\a^{k+1}_n)\left(\varphi_{r,n}(\a^{k+1}_n) \right)^{s+1} 
\right)\right|\\
&\le & C_{s,f}\,\frac{1}{n^2}\sum_{k=1}^n \left(\left(\varphi_{r,n}(\a^{k+1}_n)
\right)^{s+2}+\left(\varphi_{r,n}(\a^{k}_n) \right)^{s+2} \right)+O\left(\frac 1n\right)\\ &\le& 2C_{s,f}\, \frac 1n
\int_a^b (\varphi_{r,n})^{s+2} d  P_{r,n} +O\left(\frac 1n\right)\\ &=&O\left(\frac
1n\right). 
\end{eqnarray*}
On the other hand
\begin{eqnarray*}
\frac{1}{n} \sum_{k=1}^n\left( f(\a^k_n)\left(\varphi_{r,n}(\a^{k}_n) \right)^{s+1}+ f(\a^{k+1}_n)\left(\varphi_{r,n}(\a^{k+1}_n) \right)^{s+1} 
\right)&=& \frac 2n \int_a^b (\varphi_{r,n})^{s+1} f\,d  P_{r,n}  +O\left(\frac 1n\right)\\
&\longrightarrow& 2 \int_a^b (\varphi_{r})^{s+1} f\,d P_r\quad \mbox{as} \quad n\to \infty.
\end{eqnarray*}   
Finally
\[
\int_a^b (\varphi_{r})^{s+1} f\,d P_r= \left(\int_a^b
f^{\frac{1}{r+1}}d\lambda_{_1}\right)^{s+1}\!\!\int_a^b f^{1-\frac{s+1}{r+1}} f^{\frac{1}{r+1}}d\lambda_{_1}=
\left(\int_a^b f^{\frac{1}{r+1}}d\lambda_{_1}\right)^{s+1}\!\!\int_a^b f^{1-\frac{s}{r+1}}d\lambda_{_1}
\]
which completes the proof.$\cqfd$

\bs
\ni {\bf Remark.} In fact, we proved a slightly more precise statement, namely that, for every $r,s>0$,  
\[
\E|X-\widehat X^{\a_n}|^s= \frac{Q_{r,s}(P)}{n^s}
+O\left(\frac{1}{n^{s+1}}\right)\quad \mbox{ as } \quad n\to \infty.
\]

\section{Application to quadrature formulas for numerical
integration.}\label{Appli}
\setcounter{equation}{0}
\setcounter{Assumption}{0}
\setcounter{Theorem}{0}
\setcounter{Proposition}{0}
\setcounter{Corollary}{0}
\setcounter{Lemma}{0}
\setcounter{Definition}{0}
\setcounter{Remark}{0}

\subsection{Numerical integration on $\R^d$}
Thanks to this extension, it is now possible to take advantage of the stationarity of quadratic quantizers for a wider
class of functions. 

Let $f:\R^d\to \R$ be a twice differentiable function such that 
\[
\|D^2f(x)\|\le A(\|x\|^{d-\n}+1)
\]
for some non-negative real constants $A$ and $B$ and $\n\!\in(0,d]$. Let $X$ be a r.v. with distribution 
$P= \P_{_X}$. Assume that $X\!\in L^{2+d}(\P)$ and that for this distribution the extended exact
$L^s$-quantization rate of a sequence of optimal quadratic $n$-quantizers $(\a_n)_{n\ge1}$ holds true for every
$s\!\in (0,2+d)$.  First, let 
$\a$ be  an optimal quadratic quantizer and $\widehat X^\a$ a related Voronoi quantization of $X$. For every
couple $(p,q)\!\in(1,+\infty)^2$ of H\"older conjugate exponents, one has 
\begin{equation}\label{majorInteg}
|\E\,f(X)-\E \,f(\widehat X^{\a})| \le \frac 12 \|D^2f(\widehat
X^\a)\|_{_p}\|\,\|X-\widehat X^\a\|^2_{_{2q}}.
\end{equation}

Note that this error bound  follows from the stationarity property of a quadratic quantization that is $\E (X\,|\,
\widehat X^\a) = 
\widehat X^\a$ since 
\[
|\E\,f(X)-\E\,f(\widehat X^{\a})| = |\E\,f(X)-\E\,f(\widehat X^{\a})- \E(Df(\widehat X^\a). (X-\widehat X^\a))|\le
\frac 12\E( \|D^2f(\widehat X^\a)\|\,\| X-\widehat X^\a\|^{ 2}).
\]
 
Now,  set $p= \frac{d+2}{d-\n}$ and $q=\frac{d+2}{2+\n}$ in Equation~(\ref{majorInteg}). Then $(d-\n)p=d+2$ and
$2q<d+2$ so that, 
\[
\|D^2f(\widehat X^\a)\|_{_p}\le A'(\E\,\|\widehat X^\a\|^{d+2} +1)^{\frac 1p}\ \le A'(\E\,\| X\|^{d+2} +1)^{\frac 1p}<+\infty
\]
and
\[
 \|X-\widehat X^{\a_n}\|^2_{_{2q}}= O(n^{-\frac 2d}).
\]
Consequently the optimal rate
\[
 |\E\,f(X)-\E (f(\widehat X^{\a_n}))| = O(n^{-\frac 2d})
\]
 for numerical integration by (optimal quadratic) quantization holds for  a wide class of 
 (twice differentiable) functions whose growth at infinity can be infinitely faster than quadratic (the quadratic case  is
obtained for $q=1$ and $p=+\infty$ for functions having a bounded Hessian; then there is no need  for  the extended
quantization rate).

\subsection{Numerical integration on the Wiener space} 
Let 
\[
 W = \sum_{n\ge 1} \sqrt{\lambda_n}\,   \xi_n e^W_n
\]
the Karhunen-Lo\`eve ($K$-$L$) expansion of a standard Brownian motion $(W_t)_{t\in[0,T]}$. $(e^W_n,\lambda_n)_{n\ge 1}$ is the Karhunen-Lo\`eve
(orthonormal) eigensystem of the (nonnegative trace) covariance operator of the Brownian motion $C_{_W}(s,t)= s\wedge t$, indexed so that
$\lambda_n$ decreases (to $0$). All these quantities do have explicit expressions (see~\cite{LUS2,PAPR} among others). In that $K$-$L$ expansion, 
$(\xi_n)_{n\ge 1}$ is then a sequence of i.i.d. random ${\cal N}(0;1)$-distributed random variables.  One considers a sequence of
(scalar) product quantizations
$\widehat W^N$ of the Brownian motion
$W$
$i.e.$ 
\[
\widehat W^N = \sum_{k=1}^{m_{_N}}\sqrt{\lambda_k}\, \widehat \xi_k^{N_k} e^W_k
\]
where  for every $N\ge 1$, $N_1,\ldots,N_k\ge 1$, $1\le k\le m_{_N}$ and
$N_1\times\cdots\times N_{m_{_N}}\le N$, $\widehat \xi_k^{N_k}$ is an $L^2$-optimal $N_k$-quantization of $\xi_k$ (see~\cite{LUS1,LUS2} for more
details). Then denoting
$|\,.\,|_{L^2_{_T}}$ the quadratic norm on
$L^2([0,T],dt)$ one has  for every
$s\!\in[2,3)$,
\begin{eqnarray*}
\|W-\widehat W^N\|_{_s} &=& \E\left(|W-\widehat W|^s_{L^2_{_T}}\right)^{\frac 1s}\\
&=& \left(\E (|W-\widehat W|^2_{L^2_{_T}})^{\frac s2}\right)^{\frac 1s}\\
&=& \||W-\widehat W|^2_{L^2_{_T}}\|^{\frac 12}_{_\frac s2}\\
&=& \left\|\sum_{k=1}^{m_{_N}} \lambda_k  (\xi_k- \widehat \xi_k^{N_k})^2+\sum_{k\ge m_{_N}+1}\lambda_k \xi_k^2\right\|^{\frac 12}_{_{\frac s2}}\\
&\le & \left( \sum_{k=1}^{m_{_N}} \lambda_k  \|(\xi_k- \widehat \xi_k^{N_k})^2\|_{_\frac s2}+\sum_{k\ge m_{_N}+1}\lambda_k \|\xi_k^2\|_{_\frac s2}\right)^{\frac 12}\\
&= & \left( \sum_{k=1}^{m_{_N}} \lambda_k  \|\,\xi_k- \widehat \xi_k^{N_k}\|^2_{s}+\sum_{k\ge m_{_N}+1}\lambda_k \|\xi_k\|^2_{s}\right)^{\frac 12}
\end{eqnarray*}
where we used that $\frac s2 \ge 1$ so that $\|\,.\,\|_{_\frac s2}$ is a norm. Consequently, using that the normal distribution satisfies  
Theorem~\ref{thm2} (see Example~1), for every
$s\!\in[2,3)$ there exists a real constant $C_s$ such that
\begin{equation}\label{rsW}
\|W-\widehat W^N\|_{_2} \le \|W-\widehat W^N\|_{_s} \le C_s\left(\sum_{k=1}^{m_{_N}}\frac{\lambda_k}{N^2_k}+\sum_{k\ge m_{_N}+1}\lambda_k \right)^{\frac 12}
=C_s\|W-\widehat W^N\|_{_2}.
\end{equation}
Consequently, if furthermore $(\widehat W^N)_{N\ge 1}$ is a rate   optimal   sequence for  quadratic quantization, then 
\begin{equation}\label{rsW2}
\|W-\widehat W^N\|_{_s}\asymp \|W-\widehat W^N\|_{_2}\sim c_W (\log N)^{-\frac 12}.
\end{equation}

Note that the result clearly holds for $s\!\in(0,2]$ (simply using the monotony  of the $L^s$-norm). Furthermore, except for the final rate~$O( (\log N)^{-\frac
12})$ for which we refer to~\cite{LUS2}, the Brownian motion plays no specific r\^ole among Gaussian processes~: the inequalities in~(\ref{rsW})
extend  to the product quantization of any Gaussian process $(X_t)_{t\in[0,T]}$.

In fact, as emphasized in~\cite{LUS2},  optimal scalar product quantization is rate optimal but it is also possible to produce some rate optimal sequences
based on $d$-dimensional  block product quantizations. In fact it is shown in~\cite{LUS2} that as $d$ grows one can obtain the sharp convergence rate of
high-resolution quantization (it is even possible 
  to  choose this dimension $d=d(N)$ as a function of $N$ to achieve this sharp rate). From a numerical viewpoint   producing these
$d$-dimensional marginal blocks is more demanding than producing scalar optimal quantizations, but some issues are in progress on that topic. 

The point here is that, if one considers some $d$-dimensional marginal  blocks to produce a rate optimal sequence of quantizations of the Brownian motion, such a sequence will
satisfy~(\ref{rsW2}) for every $s\!\in(0,d+2)$ instead of $(0,3)$. This seems to be an interesting feature of $d$-dimensional block quantization.

The practical application to numerical integration on the Wiener space is formally similar to that on $\R^d$ except that one considers some
functionals $F$ on $(L^2([0,T],dt), |\,.\,|_{L^2_{_T}})$ instead of functions $f$ defined on $\R^d$. Then, one approximates $\E\,F(W)$ by
$\E\,F(\widehat W^N)$ having in mind that a closed formula exists as well for the distribution of $\widehat W^N$ ({\em the weights}, see~\cite{PAPR})
which makes the computation of
$\E\,F(\widehat W^N)$ fully tractable (with negligible cost). Then, if the Hessian $D^2 F$ of $F$ satisfies
\[
\forall\, \omega  \!\in  L^2([0,T],dt),\quad \|D^2F(\omega)\|\le A(|\omega|^r_{L^2_T}+1)
\]
for some positive real constant $A,\,r\!>\!0$  similar computations as those carried out in the $d$-dimensional case implies that 
\[
|\E\,F(W)-\E\,F(\widehat W^N)|\le \frac{C_{F,W}}{\log N}.
\] 
Without the extended quantization rate~(\ref{rsW2}) such a rate only holds for functionals $F$ with a bounded Hessian.    For further applications
(like Romberg extrapolation) and numerical experiments  devoted to path-dependent option pricing, we refer  to~\cite{PAPR}.
 
\bigskip

\bigskip
\centerline{{\bf \Large Appendix}}

In this appendix we closely follow the proof of Theorem~4 in~\cite{DEFOPA2}. This part of the proof is reproduced for the reader's convenience.

\bigskip
Set $\Delta \a_n^k=\a_n^k-\a_n^{k-1}$. After an appropriate change of variable the $L^r$-stationarity Equation~(see~\cite{DEFOPA2})
\[
r\int_{\frac{\a_n^k+\a_n^{k-1}}{2}}^{\frac{\a_n^{k+1}+\a_n^k}{2}} |a_n^k-x|^{r-1} {\rm sign}(a_n^k-x)f(x)dx=0,\qquad 1\le k\le n,
\]
 reads:
$$
\int_{-\frac{\Delta \a_n^k}{2}}^{\frac{\Delta \a_n^{k+1}}{2}}\
v^{r-1} {\rm sign}(v)\ f(\a_n^k+v)\,dv=0,\; 1\le k\le
n.
$$
Introducing $f(\a_n^k)$ we obtain:
\begin{eqnarray*}
&&\int_{-\frac{\Delta \a_n^k}{2}}^{0}\,
(-v)^{r-1}\left(f(\a_n^k+v)-f(\a_n^k)\right)dv+f(\a_n^k)
\frac{{(\Delta \a_n^k})^r}{r2^r}\\
&&\qquad=\int_{0}^{\frac{\Delta \a_n^{k+1}}{2}}\,  v^{r-1}\,
(f(\a_n^k+v)-f(\a_n^k))\,dv+f(\a_n^k)
\frac{{(\Delta \a_n^{k+1}})^r}{r2^r}.
\end{eqnarray*}

\noindent Setting $v:=\frac{\Delta\a_n^k}{2}u$ and
$v:=\frac{\Delta \a_n^{k+1}}{2}u$ respectively we have, since
$f(\a_n^k)>0$,
\begin{eqnarray*}
&&(\Delta \a_n^k)^r\left(1+r\int_0^1 u^{r-1}\frac{f(\a_n^k-\frac{\Delta
\a_n^k}{2}u)-f(\a_n^k)}{f(\a_n^k)}  \,du\right)
\\
&&\qquad =(\Delta \a_n^{k+1})^r\left(1+r\int_0^1
u^{r-1}\frac{f(\a_n^k+\frac{\Delta \a_n^{k+1}}{2}u)-
f(\a_n^k)}{f(\a_n^k)}\,du\right).
\end{eqnarray*}

\noindent Setting $\displaystyle H(x,y):=\int_0^1r
u^{r-1}\frac{f(x+u\ \frac{y-x}{2})-f(x)}{f(x)}\,du$   finally leads to
\begin{equation}\label{estim}
\left(\frac{\Delta \a_n^{k+1}}{\Delta \a_n^k}\right)^r=
\frac{1+H(\a_n^k,\a_n^{k+1})}{1+H(\a_n^k,\a_n^{k-1})}
\end{equation}

\noindent Let $[a,b]\subset (m,M)$ and let  $L^{a,b}_f$ denote the
Lipschitz coefficient of $f$ on $[a,b]$.
$$
\forall\, \xi,\,\xi' \!\in [a,b],\qquad |H(\xi,\xi')|\le
L^{a,b}_f\frac{r}{r+1} \frac{|\xi-\xi'|}{f(\xi)}\le C_{a,b,f} \,|\xi-\xi'|
$$

\noindent since
$f$ is  bounded away from $0$ on $[a,b]$. Consequently:
$$ |H(\a_n^k,\a_n^{k-1})|\le C_{a,b,f} |\a_n^k-\a_n^{k-1}|\;\mbox{ and }\;  
|H(\a_n^k,\a_n^{k+1})|\le C_{a,b,f} |\a_n^k-\a_n^{k+1}|,
$$
whenever  $\a^{k\pm 1}_n$ lie in $[a,b]$.

\noindent The $p.d.f.$ $f$ being bounded away from $0$ on $[a,b]$, 
$\max_{\{k\;/\; a\le \a^{k\pm 1}_n\le  b\}}\max(\Delta
\a_n^{k},\Delta \a_n^{k+1})$ goes to zero (see~$e.g.$Proposition~?? in~\cite{DEFOPA2}),  so we can estimate
the right hand of~(\ref{estim}):
$$
\frac{1+H(\a_n^{k},\a_n^{k-1})}{1+H(\a_n^{k},\a_n^{k+1})}
=\exp(H(\a_n^{k},\a_n^{k-1})-H(\a_n^{k},\a_n^{k+1})) 
+ O(\max(\Delta \a_n^{k},\Delta \a_n^{k+1})^2)
$$
where
\begin{equation}\label{bof0}
\!\!\!\!\!\!\!\max_{\{k\;/\; a\le \a^{k\pm 1}_n\le
b\}}\!\!\!\!|O( \max(\Delta
\a_n^{k},\Delta\a_n^{k+1})^2)| \le C\max_{\{k\;/\; a\le \a^{k\pm
1}_n\le  b\}}\!\!\!\! \max(\Delta
\a_n^{k},\Delta\a_n^{k+1})^2\stackrel{ n\to \infty}{\longrightarrow}0.
\end{equation}

\noindent Now,  we have $\displaystyle \sum_{\{k\;/\; a\le \a^{k\pm 1}_n\le  b\}}\Delta \a_n^k\le b-a$,
so that:
\begin{equation}
\displaystyle\sum_{\{k\;/\; a\le \a_n^{k\pm 1}\le  b\}}(\Delta
\a_n^k)^2\;\le \;(b-a)\hskip
-.5 cm \sup_{\{k\;/\; a\le \a_n^{k\pm 1}\le  b\}}(\Delta
\a_n^k) =(b-a)o(1).
\end{equation}
Now set
$h_\ell:=\exp(H(\a^{k+\ell}_n,\a^{k+\ell+1}_n)-H(\a^{k+\ell}_n,\a^{k+\ell-1}_n))$
and
          $\d_\ell:=O(\max(\Delta \a^{k+\ell+1}_n,\Delta \a^{k+\ell}_n)^r)$.
Note that, as long as
$\a_n^{k}$ and $\a^{\ell+p}_n\!\in[a,b]$, $\min_{1\le \ell\le p-1}
h_\ell $ is bounded away from
$0$ by a real constant $\underline{h}>0$ not depending on 
$n$  since the function $H$
is bounded away from $-\infty$ on $[a,b]^2$. Hence
$$
\left(\frac{\Delta \a^{k+p}_n}{\Delta\a_n^{k+1}}\right)^r\le
\prod_{\ell=1}^{p-1}(h_\ell+|\d_\ell|)=(\prod_{\ell=1}^{p-1}h_\ell)
\prod_{\ell=1}^{p-1}(1+\frac{\d_\ell}{h_\ell})\le(\prod_{\ell=1}^{p-1}
h_\ell)\exp(\sum_{\ell=1}^{p-1}\frac{\d_\ell}{\underline{h}})\le
(\prod_{\ell=1}^{p-1}h_\ell)\,e^{(b-a)\,o(1)}.
$$
One obtains  a lower bound the same way round. Assume that $n$ is
large enough so
that, $\max\{\d_\ell,\; a\le \a_n^{\ell\pm 1}\le  b\}\le
\underline{h}/2$.    Using that
$\ln (1-u)\ge -2u$ if $0\le u\le 1/2$,
$$
\left(\frac{\Delta \a^{k+p}_n}{\Delta\a_n^{k+1}}\right)^r\ge
(\prod_{\ell=1}^{p-1}
h_\ell)\prod_{\ell=1}^{p-1}(1-\frac{\d_\ell}{\underline{h}})
\ge (\prod_{\ell=1}^{p-1}
h_\ell)\exp(-2\sum_{\ell=1}^{p-1}\frac{\d_\ell}{\underline{h}})
\ge(\prod_{\ell=1}^{p-1}h_\ell)\,e^{(b-a)\,o(1)}.
$$
Thus it yields that for large enough  $n$ and every
$\a_{\ell+p}^n,\,\a_\ell^n\in[a,b]$:

\begin{equation}\label{estim1}
\exp(\sum_{\ell=1}^{p-1} \ln(h_\ell)-(b-a)|o(1)|)\le
\left(\frac{\Delta \a^{k+p}_n}{\Delta\a_n^{k+1}}\right)^r\le
\exp(\sum_{\ell=1}^{p-1} \ln(h_\ell)+(b-a)|o(1)|).
\end{equation}

On the other hand, using once again that  $f$ is  bounded away from
$0$ on $[a,b]$, yields
\begin{eqnarray}
\nonumber
\left| \sum_{\ell=1}^{p-1} \ln(h_\ell)\right|
          &=& \left|\sum_{\ell=1}^{p-1} r
\int_0^1\!\!\!\!u^{r-1}\frac{f(\a^{k+\ell}_n\!-\!\frac{\Delta
\a^{k+\ell}_n}{2}u)-f(\a^{k+\ell}_n+\frac{\Delta
\a^{k+\ell+1}_n}{2}u)}{f(\a^{k+\ell}_n)}du\right|\\
\label{logalpha}&\le & L^{a,b}_f \frac{r}{r+1}
\sum_{\ell=1}^{p-1}\left|\frac{\Delta \a^{k+\ell+1}_n-\Delta \a^{k+\ell}_n}{2f(\a^{k+\ell}_n)}\right|\le C(b-a).
\end{eqnarray}

\noindent  Combined with Inequality~(\ref{estim1}), this  provides
$$
\max_{\{k\;/\; a\le x_{i- 1}^{(n)}\le \a_n^{k}\le  b\}} \hskip -0.75
cm \Delta \a_n^{k}\le C\!\!\!\!\!\min_{\{k\;/\; a\le \a^{k-1}_n\le \a_n^{k}\le  b\}} \hskip -1 cm\Delta \a_n^{k}.
$$

\end{document}